\documentclass[a4paper,12pt]{elsarticle}

\usepackage[table]{xcolor}
\usepackage[english]{babel}
\usepackage{mathtools}
\usepackage{dsfont}
\usepackage{color}
\usepackage{amsfonts}
\usepackage{amsmath}
\usepackage{mathrsfs}
\usepackage{geometry}
\usepackage{multirow}
\usepackage{stmaryrd}
\usepackage{float}
\usepackage{acronym}
\definecolor{green}{RGB}{154,205,50}
\usepackage{textcomp}
\usepackage[overload]{empheq}
\usepackage{graphics,graphicx,caption}
\usepackage[capitalise]{cleveref}
\newtheorem{remark}{Remark}

\newcommand{\B}[1]{\boldsymbol{#1}}
\newcommand{\norm}[1]{|\!| #1 |\!|}
\newcommand{\jump}[1]{\mbox{$[\![ #1 ]\!]$}}

\newcommand{\FIN}{\textrm{fin}}
\newcommand{\DD}{\textsc{d}}
\newcommand{\DA}{\textsc{da}}

\renewcommand{\AA}{\textsc{a}}
\newcommand{\MAX}{\textrm{max}}
\newcommand{\OPT}{\textrm{opt}}
\newcommand{\EQ}[1]{\eqref{eq:#1}}
\newcommand{\FIG}[1]{\ref{fig:#1}}
\newcommand{\varnorm}[1]{{\left\vert\kern-0.25ex\left\vert\kern-0.25ex\left\vert #1 
    \right\vert\kern-0.25ex\right\vert\kern-0.25ex\right\vert}}

\begin{document}

\begin{frontmatter}
\title{Numerical Investigation of the Sharp-Interface Limit of the Navier-Stokes--Cahn-Hilliard Equations}

\author[TUE1]{T.H.B. Demont}
\author[TUE1]{S.K.F. Stoter}
\author[TUE1]{E.H.~van~Brummelen}

\address[TUE1]{%
  Eindhoven University of Technology, Department of Mechanical Engineering,
  P.O. Box 513, 5600 MB Eindhoven, The Netherlands}

\begin{abstract}
In this article, we study the behavior of the Abels-Garcke-Gr\"un Navier-Stokes--Cahn-Hilliard diffuse-interface model for binary-fluid flows, as the diffuse-interface thickness passes to zero. We consider this so-called sharp-interface limit in the setting of the classical oscillating-droplet problem. To provide reference limit solutions, we derive new analytical expressions for small-amplitude oscillations of a viscous droplet in a viscous ambient fluid in two dimensions. We probe the sharp-interface limit of the Navier-Stokes--Cahn-Hilliard equations by means of an adaptive finite-element method, in which the refinements are guided by an a-posteriori error-estimation procedure. The adaptive-refinement procedure enables us to consider diffuse-interface thicknesses that are significantly smaller than other relevant length scales in the droplet-oscillation problem, allowing an exploration of the asymptotic regime. For two distinct modes of oscillation, we determine the optimal scaling relation between the diffuse-interface thickness parameter and the mobility parameter. Additionally, we examine the effect of deviations from the optimal scaling of the mobility parameter on the approach of the diffuse-interface solution to the sharp-interface solution. 
\end{abstract}


\end{frontmatter}

\section{Introduction}
\label{sec:intro}

Binary-fluid flows in which the two fluid components are separated by a molecular transition layer are omnipresent in science and engineering. Examples are inkjet printing and additive manufacturing. Mathematical-physical models for binary-fluid flows generally fall under one of two categories, namely sharp-interface or diffuse-interface models. In sharp-interface models, the surface that separates the two fluid components is represented explicitly by a manifold of co-dimension one. This manifold carries kinematic and dynamic interface conditions, which act as boundary conditions on the initial boundary-value problems of the two contiguous fluid components and, in addition, determine the evolution of the manifold. Sharp-interface models are therefore of free-boundary type. In diffuse-interface models, the interface between the two fluid components is represented as a thin-but-finite transition layer, in which the two components are mixed in a proportion that varies continuously and monotonously between the two pure species across the layer. The strength of diffuse-interface models lies in their intrinsic ability to account for topological changes of the fluid-fluid interface due to coalescence or break-up of droplets or wetting, i.e.\ the propagation of the fluid-fluid front along a (possibly elastic) solid substrate~\cite{Seppecher:1996kx,Jacqmin:2000kx,Brummelen:2021aw,Yue:2011uq}.

Diffuse-interface models for two immiscible incompressible fluid species are generally described by the Navier-Stokes--Cahn-Hilliard (NSCH) equations. The NSCH equations represent a class of models, of which various renditions have been proposed over the last half century: by Hohenberg and Halperin in the late 1970s \cite{Hohenberg:1977hh}, by Lowengrub and Truskinovsky in the late 1990s \cite{Lowengrub:1998uq}, by Shokrpour et al. in 2018 \cite{Simsek:2018gb} and by Abels, Garcke and Gr\"{u}n in 2012 \cite{Abels:2012vn}.
In this article, we focus on the latter model, in view of its thermodynamic consistency and its consistent reduction to the underlying single-fluid Navier--Stokes equations in the pure species setting.

NSCH models invariably contain three parameters related to the diffuse interface,  viz.\ an interface-thickness parameter, $\varepsilon$, a mobility parameter, $m$, and a surface-tension parameter, $\sigma$. The interface-thickness parameter represents the transverse length scale of the transition layer between the two fluid components, and the transition layer collapses (specifically, is supposed to collapse) onto a manifold of co-dimension one in the so-called sharp-interface limit $\varepsilon\to+0$. The mobility parameter is responsible for the rate at which phase-diffusion occurs in the vicinity of the diffuse interface. In the phase-separated regime in which the NSCH equations are typically applied as a binary-fluid model, the mobility parameter is responsible for the rate at which the interface equilibrates. In the mixture regime, it governs the dynamics of the Ostwald-ripening effect. The surface-tension parameter controls the excess free energy $\sigma_{\DA}$ of the diffuse-interface 
according to~$2\sqrt{2}\sigma=3\sigma_{\DA}$. It is to be noted that for the NSCH equations this proportionality holds independent of $\varepsilon$, as opposed to the Navier-Stokes--Korteweg equations.

Contemporary understanding of the sharp-interface limit of the NSCH equations is incomplete. An overview of known results and open questions is provided in~\cite[\S{}4.3]{Abels:2018ly}. One prominent open question pertains to the appropriate scaling of the mobility parameter in relation to the interface-thickness parameter, in the sharp-interface limit. The limit solution of the NSCH equations depends on the scaling $m:=m_{\varepsilon}$. Abels and Garcke~\cite[\S{}4.1]{Abels:2018ly} establish that, 
if $m_{\varepsilon}=\mathcal{O}(1)$ as $\varepsilon \rightarrow +0$, their NSCH model converges to the nonclassical sharp-interface Navier--Stokes/Mullins--Sekerka model; see also~\cite[\S{}4]{Jacqmin:1999fk}. If, on the other hand, $m_{\varepsilon}$ vanishes suitably as $\varepsilon\to+0$, the classical sharp-interface binary-fluid model is obtained, where the interface is transported by the fluid velocity \cite{Lowengrub:1998uq,Abels:2018ly,Yue:2010hq}. However, the decay of the mobility cannot be too fast: if $m_{\varepsilon}=o(\varepsilon^3)$ as $\varepsilon\to+0$, the resulting limit solution of the NSCH model generally violates the Young--Laplace condition on the pressure jump across the interface~\cite{Abels:2014ca}. These results suggest that $m_{\varepsilon}\propto{}\varepsilon^{a}$ with $0<a\leq{}3$ as $\varepsilon\to+0$ represents a necessary and sufficient condition to achieve a classical sharp-interface solution. Still, the details of the approach of the diffuse-interface solution to the sharp-interface limit solution for these various admissible scalings of the mobility are not currently known, and different scalings have been proposed in the literature,
in particular in the context of numerical simulation approaches. In~\cite{Demont:2022dk}, the scaling
$m_{\varepsilon} \propto \varepsilon^3$ is
 considered, based on the argument that this proportionality fixes the diffusive time scale and, thus, the equilibration rate of the diffuse interface. This cubic scaling of the mobility with respect to the interface thickness (in terms of their usual dimensional forms) is also propounded in~\cite{Khatavkar:2006gf}, supported by numerical investigations. Based on partial matched-asymptotic-analysis arguments, Ref.~\cite{Magaletti:2013vn} concludes that $m \propto \varepsilon^2$ is the appropriate scaling. However, because the matching procedure in this reference is incomplete, it is unclear whether this scaling relation in fact represents a necessary or sufficient condition. On the basis of a consideration of curvature-induced expansion/contraction modes at the diffuse interface, it is argued in~\cite{Jacqmin:1999fk} that $m_{\varepsilon} \propto \varepsilon^a$ with $1\leq{}a<2$. It is to be noted that the aforementioned scalings of the mobility pertain to situations without moving contact lines and topological changes; see, e.g.~\cite{Yue:2010hq}.
 
 In this article, we address the open questions of the optimal scaling of the mobility and the approach to the sharp-interface limit solution by computational investigation of the behavior of the Abels--Garcke--Gr\"un Navier-Stokes--Cahn-Hilliard model for different interface dynamics and different mobility parameters as it limits toward a sharp-interface description of a two-dimensional oscillating droplet. To enable an exploration of the asymptotic regime, we apply an adaptive finite-element method, in which the adaptivity is guided by an a-posteriori error estimate; see~\cite{Brummelen:2021aw,Demont:2022dk} for details. 

We conduct our analysis of the sharp-interface limit of the NSCH equations in the context of the prototypical oscillating-droplet problem, in two dimensions. Despite the fact that the oscillating-droplet problem is classical, it appears that the two-dimensional setting has not been extensively investigated, and that solutions of the two-dimensional problem have not been reported in the literature. The investigation of the oscillating-droplet problem dates back to Rayleigh, who presented the well-known frequency of oscillation of an inviscid droplet {\em in vacuo\/} in 1879~\cite{Strutt:1879xi}. This result was extended by Lamb in the 1930s to include the effect of an inviscid ambient fluid~\cite{Lamb:1932ss}. In 1960, Reid generalized the theory of oscillating droplets {\em in vacuo\/} by including the effect of viscosity~\cite{reid1960oscillations}. A complete theory, comprising solutions for small oscillations of a viscous droplet in a viscous ambient fluid,  was then finally presented by Miller and Scriven in 1968~\cite{Miller:1968cu}. The aforementioned references however exclusively consider the three-dimensional case and the results, especially those for the viscous solutions, do not trivially extend to the two-dimensional case. Clearly, the three-dimensional case is the practically relevant one, but the two-dimensional case has {\em raison d'\^etre\/} independently as a means of verification for mathematical models and numerical methods. Our analysis of the sharp-interface limit of the NSCH equations requires access to closed form solutions of the sharp-interface model, on the one hand to provide initial and boundary data for the NSCH equations, and on the other hand to systematically determine the deviation of the diffuse-interface solution relative to the sharp-interface solution. A secondary objective of this work is therefore to establish closed-form expressions for small-amplitude oscillations of a viscous droplet in a viscous ambient fluid. Our derivation follows that of Miller and Scriven, but we deviate from their derivation by a more complete elaboration of intermediate steps and assumptions and, in particular, an explicit accounting of the complex-valued nature of the different fields, and by presenting closed-form expressions of the final results.

The remainder of this article is structured as follows. In \cref{sec:prob_form}, we lay out the Abels--Garcke--Gr\"un Navier-Stokes--Cahn-Hilliard model equations, and the coupled Navier--Stokes free-boundary problem that they should reduce to in the sharp-interface limit. In \cref{sec:deriv}, we derive a closed form expression for the sharp-interface model corresponding to  small-amplitude oscillations of a viscous droplet in a viscous ambient fluid in two dimensions. We make use of these expressions in \cref{sec:num_exp}, where we study the approach of the NSCH solution to the sharp-interface solution in the limit $\varepsilon\to+0$, by means of systematic numerical experiments.

%===============================================================================================================%
%===============================================================================================================%
%===============================================================================================================%

\section{Governing equations}
\label{sec:prob_form}

\begin{figure}[!b]
\begin{center}
\includegraphics[width=0.65\textwidth]{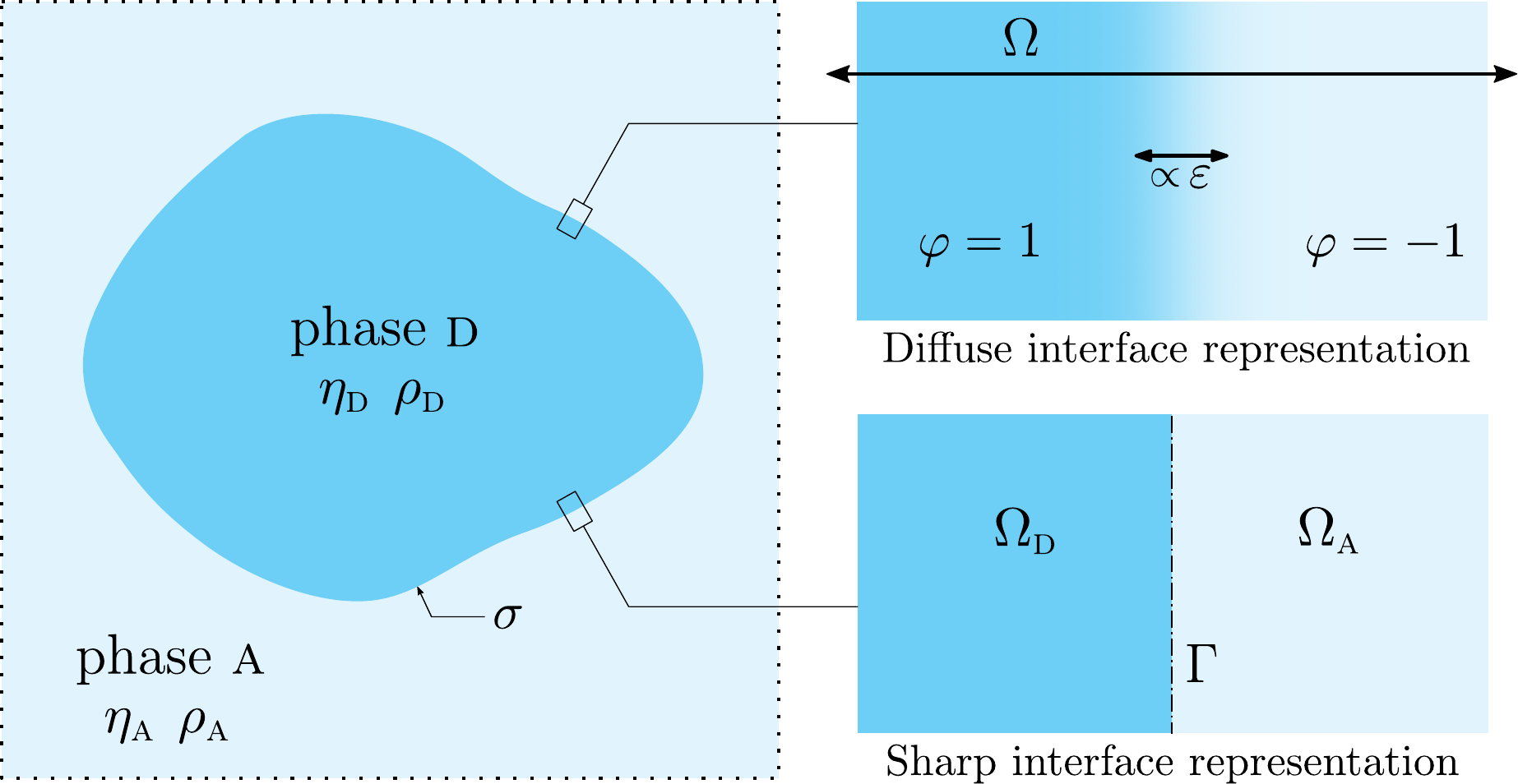}
\end{center}
\caption{Schematic of the physical setting of a submerged immiscible fluid, modeled with either a diffuse interface representation or a sharp-interface representation.}
\label{fig:droplet_model}
\end{figure}

We consider a binary-fluid system, where both fluids are modeled as being incompressible, isothermal, immiscible, and Newtonian with finite viscosity. In accordance with the later focus on a submerged droplet, we denote one of the fluids by $\DD$, for droplet, and the other by $\AA$, for ambient. Various modeling frameworks for describing the fluid motion exist. These make use of either a diffuse interface representation or a sharp-interface representation. \Cref{fig:droplet_model} illustrates the different relevant domains and material parameters. We consider in this work the incompressible NSCH model --- specifically, the model developed by Abels, Garcke, and Gr\"{u}n in \cite{Abels:2012vn} --- in order to describe the binary fluid dynamics. Motivation for this choice lies in its thermodynamic consistency, incompressibility, and consistent reduction to the underlying single-fluid Navier--Stokes equations in the pure species setting. Recently, the well-posedness of the Abels--Garcke-Gr\"{u}n model in various settings has been 
shown~\cite{Giorgini:2021bg, Abels:2013bd, Abels:2013pa}.

\subsection{Diffuse-interface representation}
\label{sec:dif_int_model}

In diffuse-interface models, the two immiscible fluids are separated by a layer of finite thickness constituted by a mixture of both fluids, reflecting a gradual transition between fluid $\DD$ and fluid $\AA$. We consider an open time interval $(0, t_\FIN) \subseteq \mathds{R}_{>0}$ and a spatial domain corresponding to a simply connected time-independent subset $\Omega \subseteq \mathds{R}^d$ ($d=2,3$). We make use of a Navier-Stokes--Cahn-Hilliard type model that describes the evolution of a so-called order parameter $\varphi \in [-1,1]$ representing pure species $\DD$ and $\AA$ when $\phi=1$ and $\phi=-1$, respectively, and a mixture of both when $\varphi \in (-1,1)$, in addition to the velocity and pressure of the mixture. The NSCH model as presented by Abels, Garcke, and Gr\"{u}n 
is given by~\cite{Abels:2012vn}:
\begin{subequations} \label{eq:strong}
\begin{empheq}[right={\empheqrbrace \textrm{ in } \Omega}]{align}
\partial_t \left( \rho \B{u} \right) + \nabla \cdot \left( \rho \B{u} \otimes \B{u} \right) + \nabla \cdot \left( \B{u} \otimes \B{J} \right) + \nabla p - \nabla \cdot \B{\tau} - \nabla \cdot \B{\zeta} & = 0,\\
\nabla \cdot \B{u} & = 0,\\
\partial_t \varphi + \nabla \cdot \left( \varphi \B{u} \right) - \nabla \cdot \left( m \nabla \mu \right) & = 0,\\
\mu + \sigma \varepsilon \Delta \varphi - \frac{ \sigma } { \varepsilon } \Psi' & = 0,
\end{empheq}
\end{subequations}
where the volume-averaged velocity $\B{u}$, the pressure $p$, the order parameter $\varphi$ and the chemical potential $\mu$ are the unknown fields. The closure relations for the relative mass flux $\B{J}$, the viscous stress $\B{\tau}$, the capillary stress $\B{\zeta}$ and the mixture energy density $\Psi$ are given as:
\begin{subequations}\label{eq:closureeqs}
\begin{alignat}{3}
& \B{J} \coloneqq m \frac{ \rho_{\AA} - \rho_{\DD} } { 2 } \nabla \mu  \,,\\
& \B{\tau} \coloneqq \eta ( \nabla \B{u} + ( \nabla \B{u} )^T ) \,, \\
& \B{\zeta} \coloneqq - \sigma \varepsilon \nabla \varphi \otimes \nabla \varphi + \B{I} \left( \frac{ \sigma \varepsilon } { 2 } | \nabla \varphi |^2 + \frac{ \sigma } { \varepsilon } \Psi \right) \,, \\
& \Psi \left( \varphi \right) \coloneqq \frac{1}{4} \left( \varphi^2 -1 \right)^2 \,.
\end{alignat}
\end{subequations}
The remaining parameters are material and model parameters. The model parameters are the mobility parameter $m>0$ and the interface thickness parameter $\varepsilon>0$, which affect the time and length scale of the diffuse interface, respectively. The material parameters are $\sigma$, a rescaling of the droplet-ambient surface tension $\sigma_\DA$ according to $2\sqrt{2}\sigma=3\sigma_\DA$, the mixture density $\rho$, and the mixture viscosity $\eta$. The mixture density and viscosity generally depend on~$\varphi$.
To ensure existence of a solution to the system of equations, we must allow $\varphi$ to take on values outside of $[-1,1]$~\cite{Grun:2016gi}. We include a density extension that ensures positive densities even for the nonphysical scenario $\varphi \notin [-1,1]$~\cite{Bonart:2019re}:
\begin{equation}\label{eq:densityextension}
\rho(\varphi) = \left\{
\begin{tabular}{ll}
$\frac{ 1 } { 4 } \rho_{\AA},$ & $\varphi \leq - 1 - 2 \lambda \,,$\\
$\frac{ 1 } { 4 } \rho_{\AA} + \frac{ 1 } { 4 } \rho_{\AA} \lambda^{-2} \left( 1 + 2 \lambda + \varphi \right)^2,$ & $\varphi \in ( - 1 - 2 \lambda , - 1 - \lambda) \,,$\\
$\frac{ 1 + \varphi } { 2 } \rho_{\DD} + \frac{ 1 - \varphi } { 2 } \rho_{\AA},$ & $\varphi \in [ - 1 - \lambda, 1 + \lambda ] \,,$\\
$\rho_{\DD} + \frac{ 3 } { 4 } \rho_{\AA} - \frac{ 1 } { 4 } \rho_{\AA} \lambda^{-2} \left( 1 + 2 \lambda - \varphi \right)^2,$ & $\varphi \in ( 1 + \lambda, 1 + 2 \lambda ) \,,$\\
$\rho_{\DD} + \frac{ 3 } { 4 } \rho_{\AA},$ & $\varphi \geq 1 + 2 \lambda \,,$
\end{tabular}
\right.
\end{equation}
where $\lambda = \rho_{\AA} / \left( \rho_{\DD} - \rho_{\AA} \right)$. For the viscosity interpolation, we apply the Arrhenius mixture-viscosity model~\cite{Arrhenius:1887xr}:
\begin{equation}\label{eq:arrhenius}
\log \eta( \varphi ) = \frac{ \left( 1 + \varphi \right) \Lambda \log \eta_{\DD} + \left( 1 - \varphi \right) \log \eta_{\AA} } { \left( 1 + \varphi \right) \Lambda + \left( 1 - \varphi \right) } \,,
\end{equation}
where $\Lambda = \frac{ \rho_{\DD} M_{\AA} } { \rho_{\AA} M_{\DD} }$ is the intrinsic volume ratio between the two fluids (with $M_{\AA}$ and $M_{\DD}$ their respective molar masses).

\begin{remark}To eliminate the Ostwald-ripening effect in the pure species, a degenerate dependence of the mobility on the phase field can be introduced, according to $m(\varphi) \geq 0$ with inequality if and only if $|\varphi|<1$. 
However, as a degenerate mobility introduces complications with regard to numerical-approximation procedures~\cite{Barrett:1999nx}, we opt for a constant mobility parameter.
\end{remark}

\begin{remark}
In this work, we use a volume-fraction-based Arrhenius relation, i.e.\ $\Lambda=1$. Because the denominator in~\cref{eq:arrhenius} then reduces to a non-zero constant, this choice eliminates singularities, and the mixture viscosity is bounded away from zero in a finite interval including $[-1,1]$. See \emph{Remark 2} in \cite{Brummelen:2021aw} for further details.
\end{remark}

\subsection{Sharp-interface limit}
As $\varepsilon \rightarrow +0$, the width of the diffuse interface in the NSCH model reduces to zero. As pointed out in the introduction, the particular model that arises in this limit depends on the scaling relation of the mobility~$m$. If the mobility also tends to zero appropriately, then the following classical sharp-interface model is obtained:
\begin{subequations}
\label{eq:sharp_strong}
\begin{alignat}{3}
\rho_{i}\partial_t \B{u}_i + \rho_{i} \left( \B{u}_i \cdot \nabla \right) \B{u}_i - \eta_{i} \Delta \B{u}_i + \nabla p_i & = 0 && \textrm{ in } \Omega_{i} = \Omega_{i}(t) \,, \label{eq:sharp_strong_a}\\
\nabla \cdot \B{u}_i & = 0 && \textrm{ in } \Omega_{i} \,, \label{eq:sharp_strong_b}\\
\B{u}_{i}\cdot\B{n} & = \mathcal{V}\quad && \textrm{ on } \Gamma = \Gamma(t) \,, \label{eq:IF_normal}\\
\jump{\B{u}\cdot \B{t}_j} & = 0 \quad && \textrm{ on } \Gamma \text{, for }j=1,\cdots,d-1 \,,\\
\jump{ -( \nabla \B{u} + ( \nabla \B{u} )^T ) \B{n} + p \B{n} } & = \sigma_{\DD\AA} \kappa \B{n} \quad&& \textrm{ on } \Gamma  \,, \label{eq:sharp_strong_e}
\end{alignat}
\end{subequations}
for $i \in \{\DD,\AA\}$, and with $\B{n}$ the unit normal vector on $\Gamma$ external to~$\Omega_{\DD}$, $\mathcal{V}$ the interface normal velocity, $\kappa$ the (additive) curvature of the interface, and $\jump{\cdot}$ the interface jump operator $\jump{g} = g|_{\AA} - g|_{\DD}$. We adhere to the convention that curvature is negative if the center of the osculating circle in the normal plane is located in the droplet domain.

As opposed to the NSCH model, the sharp-interface model~\eqref{eq:sharp_strong}
represents a set of equations for each fluid species separately, complemented by appropriate coupling conditions at the interface. The sharp-interface model 
represents a free-boundary problem. The domains on which the various fields are defined evolve in time, as reflected by the time-dependence of $\Omega_i(t)$ and $\Gamma(t)$. There is an intrinsic coupling between the velocity field and the evolution of $\Omega_i$ and $\Gamma$, according to~\eqref{eq:IF_normal}. Because this equation holds on both sides of the interface and~$\mathcal{V}$ is single-valued, \cref{eq:IF_normal} implies~$\jump{\B{u} \cdot \B{n} } = 0$.

%===============================================================================================================%
%===============================================================================================================%
%===============================================================================================================%

\section{Response of an oscillating droplet}
\label{sec:deriv}
With the objective of providing a reference solution for the sharp-interface limit of the NSCH equations, we now consider solutions of the free-boundary  problem~\eqref{eq:sharp_strong} corresponding to small perturbations of a circular droplet set in an ambient fluid, in two dimensions. Denoting by~$R_0$ the radius of the droplet, one can verify that
\begin{subequations}
\label{eq:gensol}
\begin{alignat}{2}
\Omega_{\DD,0}&=\{\boldsymbol{x}\in\mathbb{R}^2:|\boldsymbol{x}|<R_0\}
\label{eq:gensol_a}
&\qquad
\Omega_{\AA,0}&=\mathbb{R}^2\setminus\overline{\Omega_{\DD,0}}\,,
\\
\boldsymbol{u}_{\DD,0}&=0&\qquad\boldsymbol{u}_{\AA,0}&=0\,,
\\
p_{\DD,0}&=\sigma_{\DD\AA}/R_0\
&\qquad
p_{\AA,0}&=0\,,
\end{alignat}
\end{subequations}
represents a stationary solution to~\eqref{eq:sharp_strong}. We will use~\eqref{eq:gensol} as a generating solution. We consider perturbations of the solution~\eqref{eq:gensol} that are suitably bounded and vanish toward infinity. 
The latter condition can be expressed as:
\begin{equation}
\label{eq:BCs}
\lim\limits_{|\B{x}|\rightarrow\infty} \big(\B{u},p\big)(\B{x},t) =  0 \,.
\end{equation}
Our derivation of the natural response of such a droplet follows that of Miller and Scriven~\cite{Miller:1968cu}, except that we provide a more complete elaboration of intermediate steps and assumptions and, in particular, an explicit accounting of the complex-valued nature of the different fields. The approach essentially comprises four steps. First, we linearize the governing equations around the generating solution~\eqref{eq:gensol}, perturbed by a small deformation of the interface. Second, we derive the general solutions corresponding to the natural response in both domains separately. In the third step, we incorporate the interface coupling conditions by constraining the free parameters in the general solutions. Finally, the characteristic temporal response (frequency of oscillation and rate of damping) of the assumed interface displacement, as well as the corresponding shapes, follow from a solution-existence condition.

%===============================================================================================================%

\subsection{Formal linearization}
\label{ssec:lineariztion}
We consider small-amplitude perturbations of the interface that are sinusoidal along the circumference of the droplet. To facilitate the presentation, we introduce polar coordinates $r\in\mathbb{R}_{\geq{}0}$ and $\theta\in[0,2\pi)$ and the coordinate transformation $\B{x}=(x_1,x_2)=r(\cos\theta,\sin\theta)$. 
We regard perturbations of the interface $\Gamma_0=\partial\Omega_{\DD,0}$ corresponding to the following parametrization:
\begin{equation}
\label{eq:Gamma_delta}
\Gamma_{\delta}(t)=\big\{\boldsymbol{x}\in\mathbb{R}^2:
\boldsymbol{x}=R_{\delta}(\theta,t)\,(\cos(\theta),\sin(\theta)),\theta\in[0,2\pi)\big\}  \,,
\end{equation}
where
\begin{equation}
\label{eq:R_delta}
\begin{aligned}
R_{\delta} ( \theta , t ) & = R_0 + R_0\, \delta\left( \beta \cos ( k \theta ) + \sqrt{1-\beta^2 } \, \sin  ( k \theta ) \right) \cos(\nu t) e^{- \alpha t} \\
& = \Re \left( R_0 + R_0\, \delta \left( \beta \cos ( k \theta ) + \sqrt{1-\beta^2 } \, \sin ( k \theta ) \right) e^{- \gamma t} \right) \,.
\end{aligned}
\end{equation}
The interface configuration~\eqref{eq:Gamma_delta}-\eqref{eq:R_delta} represents a damped oscillation of the droplet with a mode-shape described by the mode number $k \in \mathbb{N}$, with angular orientation dependent on $0\leq \beta\leq 1$, and with an amplitude described by $\delta\ll{}1$ as the fraction of the droplet radius $R_0$. Our interest is restricted to droplet configurations for which $\operatorname{meas}(\Omega_{\DD})=\operatorname{meas}(\Omega_{\DD,0})+\mathcal{O}(\delta^2)$ and the barycenter of $\Omega_{\DD}$ is located at the origin. This implies that $k \in \mathbb{N}_{\geq2}$. The damping rate $\alpha \geq 0$ and the frequency of oscillation $\nu > 0$ implicitly depend on the mode number and will follow from the subsequent analysis. The second expression in~\eqref{eq:R_delta} provides a representation of $R_{\delta}$ as the real part of a complex-valued function,
with $\gamma \coloneqq \alpha - i \nu$. This form enables us to condense some of the expressions that appear in the sequel. 

In conjunction with the interface configuration~\eqref{eq:Gamma_delta}-\eqref{eq:R_delta}, we consider linear asymptotic solutions of the sharp-interface problem of the form 
\begin{equation}
\label{eq:pertsol}
\big(\B{u}_{i},p_i,\B{n}, \B{t},\kappa,\mathcal{V})=
\big(\B{u}_{i},p_i,\B{n}, \B{t},\kappa,\mathcal{V})_0
+
\delta\big(\B{u}_{i},p_i,\B{n}, \B{t},\kappa,\mathcal{V})_1\,,
\end{equation}
i.e.\ functions conforming to~\eqref{eq:pertsol} that satisfy~\eqref{eq:sharp_strong} modulo terms of~$o(\delta)$ as $\delta\to{}0$.
Substituting~\eqref{eq:pertsol} into the sharp-interface equations~\eqref{eq:sharp_strong}, collecting terms of distinct orders in~$\delta$, and noting that all terms of~$\mathcal{O}(1)$ vanish on account of the fact that the first term in~\eqref{eq:pertsol} represents a solution to~\eqref{eq:sharp_strong}, we obtain the following infinitesimal conditions on the second term in~\eqref{eq:pertsol}:
\begin{subequations}
\label{eq:sharp_1st}
\begin{alignat}{3}
\rho\partial_t \B{u}_{i,1} - \eta \Delta \B{u}_{i,1} + \nabla p_{i,1} & = 0 && \textrm{ in } \Omega_{i,0} \,, \label{eq:sharp_1st_a}\\
\nabla \cdot \B{u}_{i,1} & = 0 && \textrm{ in } \Omega_{i,0} \,, \label{eq:sharp_1st_b}\\
 \B{u}_{i,1} \cdot \B{n}_0 &= \mathcal{V}_1 && \textrm{ on } \Gamma_0 \,, \label{eq:sharp_1st_c}\\
\jump{\B{u}_{1} \cdot \B{t}_0} & = 0 && \textrm{ on } \Gamma_0 \,, \label{eq:sharp_1st_d}\\[-0.15cm]
\jump{ - \eta \left(\nabla \B{u}_{1} + \left(\nabla \B{u}_{1} \right)^T \right) \B{n}_0 + p_{1} \B{n}_0 } & = \sigma_{\DD\AA} \kappa_1 \B{n}_0 \qquad && \textrm{ on } \Gamma_0 \,, \label{eq:sharp_1st_e}
\end{alignat}
\end{subequations}
for $i \in \{\DD,\AA\}$. The non-linear advective term has dropped since the only non-linear first-order perturbation terms are cross terms between $\B{u}_{i,1}$ and $\B{u}_{i,0} = 0$. Similarly, only $\B{n}_0$ appears in the first-order conditions~\eqref{eq:sharp_1st}, because its multiplication with $\B{u}_{i,1}$ is a second order term, its dot-product with $\B{u}_{i,0}$ vanishes, and $\jump{p_0 \B{n}_1} = \sigma_{\DD\AA} R^{-1}_0 \B{n}_1$ cancels with the right-hand-side $\sigma_{\DD\AA} \kappa_0 \B{n}_1$. 

\begin{remark}
The first-order conditions are set on the stationary generating domains, $\Omega_{i,0}$, and the generating interface, $\Gamma_0$. This is a universal characteristic of linearizations of free-boundary problems. The perturbation of the interface according to~\eqref{eq:R_delta} appears implicitly in~\eqref{eq:sharp_1st} in the interface-velocity perturbation, $\mathcal{V}_1$, and the curvature perturbation, $\kappa_1$.
\end{remark}

The remainder of this section is devoted to finding general solutions to \eqref{eq:sharp_1st} for different perturbation wave-numbers~$k$. For purposes of readability, henceforth we suppress the subscripts corresponding to the order of perturbation.

%===============================================================================================================%

\subsection{First-order solutions in the droplet and ambient domain}
Next, we derive the general first-order solutions $(\B{u},p)_{i,1}$ in accordance with the differential equations~\eqref{eq:sharp_1st_a}--\eqref{eq:sharp_1st_b} for both the droplet $\DD$ and ambient $\AA$ domains. We proceed by deriving the vorticity equation corresponding to~\eqref{eq:sharp_1st_a}, which we solve by means of separation of variables. The first-order velocity field is subsequently retrieved from the vorticity solutions. The corresponding first-order pressure fields are derived as those that yield balance of linear momentum.

\subsubsection{Vorticity solution}
The pressure may be eliminated from the governing equations by taking the curl of~\eqref{eq:sharp_1st_a} and by using the identity $\nabla \times \nabla(\cdot) = 0$. Introducing the vorticity~$\omega=\nabla \times \B{u}$, we infer from~\eqref{eq:sharp_1st_a} that
\begin{align} \label{eq:NS_curl}
\rho \partial_t \left( \nabla \times \B{u} \right) - \eta \Delta \left( \nabla \times \B{u} \right) + \nabla \times\nabla{}p= \rho \partial_t \omega - \eta \Delta \omega = 0\,.
\end{align}
Let us note that in a two-dimensional setting, vorticity can be represented as a scalar-valued field. The evolution equation for this scalar vorticity field can be recognized as a diffusion equation.

To determine the general solution to~\eqref{eq:NS_curl}, we assume the following separation of variables form:
\begin{equation}
\label{eq:sepvar}
\omega ( r,\theta , t ) \coloneqq \Psi(r)\, \Theta(\theta)\, T(t).
\end{equation}
Substitution in the polar coordinate representation of the diffusion equation gives
\begin{align}
    \Psi(r)\, \Theta(\theta)\, T'(t) = \frac{\eta}{\rho r^2}\Psi(r)\, \Theta''(\theta)\, T(t)+ \frac{\eta}{\rho r}\Psi'(r)\, \Theta(\theta)\, T(t) + \frac{\eta}{\rho}\Psi''(r)\, \Theta(\theta)\, T(t) \,, \label{polarheat}
\end{align}
where primes denote differentiation. The usual separation of variables argument leads to
\begin{subequations}
\label{eq:sepvararg}
\begin{alignat}{2}
  & T'(t) = -\frac{\eta}{ \rho} m^2\, T(t) \quad\,\,\,  m \in \mathbb{C}\,, \label{eq:sepvararga}\\[-0.05cm]
  &  \Theta''(\theta) = -n^2 \, \Theta(\theta)  \qquad n \in  \mathbb{C} \,,  \label{eq:sepvarargb} \\[0.1cm]
  &  r^2 \Psi''(r) + r \Psi'(r) + \left(m^2 r^2 - n^2 \right) \Psi(r)  = 0 \,, \label{eq:BesselODE}
\end{alignat}
\end{subequations}
with $\mathbb{C}$ the set of complex numbers. 

The general solutions of~\eqref{eq:sepvararga} and~\eqref{eq:sepvarargb} consist of complex-valued exponential functions, according to
\begin{subequations}
\label{eq:ansatz}
\begin{alignat}{2}
&T (t) = c\, e^{- \frac{\eta}{ \rho} m^2 t}\,, \label{eq:tansatz} \\
&\Theta (\theta) = c_1\, e^{ i n t} + c_2 e^{ -i n t}\,.\label{eq:thetaexp}
\end{alignat}
\end{subequations}
From the periodicity of the droplet perturbations in the angular dependence, conforming to~\eqref{eq:Gamma_delta}, we infer that $n\in\mathbb{Z}_{\geq{}0}$. The arbitrary constants~$c_1$ and~$c_2$ can then be selected such that~\eqref{eq:thetaexp} reduces to the sum of two real-valued trigonometric functions:
\begin{equation}
\label{eq:thetaansatz}
\Theta (\theta) = C\cos( n \theta) +  D\sin( n \theta)\qquad (n \in \mathbb{Z}_{\geq0}) \,, 
\end{equation}
where $C,D \in \mathbb{R}$ are coefficients that determine the angular orientation of the solution. Regarding~\eqref{eq:BesselODE}, we note that this equation corresponds to Bessel's equation with a complex-valued scaling~$m\in\mathbb{C}$. Solutions of~\eqref{eq:BesselODE} therefore consist of extensions of Bessel functions to the complex plane. Such extensions of Bessel functions are well defined, by virtue of the fact that Bessel functions are analytic functions on~$\mathbb{R}$ and can hence be extended to analytic functions on~$\mathbb{C}$
via their power-series expansion. The general solution of~\eqref{eq:BesselODE} consists of a linear combination of two Bessel functions of order~$n\in\mathbb{Z}_{\geq0}$. For reasons that will become clear once we consider the boundary conditions, we choose to work with a Bessel function of the first kind, $J_n$, and a Hankel function of the second kind, $H^{(2)}_n$:
\begin{equation}
\label{eq:Psisol}
\Psi(r) = A m^2 J_n ( m r ) + B m^2 H^{(2)}_n ( m r ) \,, 
\end{equation}
with $A,B \in \mathbb{C}$. Substitution of~\eqref{eq:tansatz}, \eqref{eq:thetaansatz} and~\eqref{eq:Psisol} into~\eqref{eq:sepvar} gives the general rotationally periodic solution of the vorticity equation~\eqref{eq:NS_curl}:
\begin{equation}
\label{eq:genomega}
\omega(r,\theta,t) = \Big[ A m^2 J_n ( m r ) + B m^2 H^{(2)}_n ( m r ) \Big] \Big[ C\cos( n \theta) +  D\sin( n \theta) \Big]\, e^{- \frac{\eta}{ \rho} m^2 t} \,,
\end{equation}
for arbitrary $A,B\in\mathbb{C}$, $C,D\in\mathbb{R}$ and $m\in\mathbb{C}$.

%===============================================================================================================%

\subsubsection{Velocity solutions }
To obtain the velocity fields from the general vorticity solution~\eqref{eq:genomega}, we introduce a stream function $\chi$ according to:
\begin{equation}
    \label{eq:4}
    \Delta \chi = - \omega \,.
\end{equation}
The velocity can be retrieved from this stream function as
\begin{equation}
\label{eq:uchi}
   \B{u} = \frac{1}{r} \frac{\partial}{\partial\theta} \chi\, \B{e}_r - \frac{\partial}{\partial r} \chi \, \B{e}_\theta\,. 
\end{equation}
Based on the expression for $\omega$ in \eqref{eq:genomega}, we anticipate that a particular solution to \eqref{eq:4} is of the form:
\begin{equation}
\label{eq:4a}
 \chi (r,\theta,t) = \Upsilon ( r ) \, \big[ C\cos( n \theta) +  D\sin( n \theta) \big]\, e^{- \frac{\eta}{ \rho} m^2 t} \,.
\end{equation}
Substitution of~\eqref{eq:4a} into~\eqref{eq:4} leads to the following ordinary differential equation for~$\Upsilon$:
\begin{equation}
\label{eq:5}
\Upsilon'' ( r ) + \frac{1}{r} \Upsilon' ( r ) - \frac{n^2}{r^2} \Upsilon ( r ) = - A m^2 J_n ( m r ) - B m^2 H^{(2)}_n ( m r ) \,.
\end{equation}
The general solution to the nonhomogeneous ordinary differential equation~\eqref{eq:5} is given by:
\begin{equation}
\label{eq:Upsilon}
\Upsilon(r)   = E r^{n} + F r^{-n} + A J_n ( m r) + B H^{(2)}_{n} ( m r ) \,.
\end{equation}
The first and second term in~\eqref{eq:Upsilon} constitute the homogeneous part of the solution. The third and fourth term represent the particular part.
From \eqref{eq:uchi} it then follows that the $r$ and $\theta$ components of the velocity field are given by:
\begin{subequations}
\label{eq:usol}
\begin{alignat}{2}
&u_r = \frac{1}{r} \frac{\partial}{\partial\theta} \chi = \frac{n}{r} \Upsilon ( r ) \, \big[ D \cos( n \theta) -  C \sin( n \theta) \big]\, e^{- \frac{\eta}{ \rho} m^2 t} \,,\\
&u_\theta =  - \frac{\partial}{\partial r} \chi = - \Upsilon' ( r ) \, \big[ C \cos( n \theta) +  D \sin( n \theta) \big]\, e^{- \frac{\eta}{ \rho} m^2 t} \,.
\end{alignat}
\end{subequations}

The velocity solutions~\eqref{eq:usol} are not generally bounded in the limits
$r\to{}0$, $r\to\infty$ and $t\to\infty$. Auxiliary conditions must be imposed on the coefficients in~\eqref{eq:usol} to extract general solutions that are bounded in the droplet domain $\Omega_{\DD,0}$ (resp. ambient domain $\mathbb{R}^2\setminus\overline{\Omega_{\DD,0}}$) as $r\to{}0$ (resp.~$r\to\infty$) and as $t\to\infty$. To ensure boundedness of the solutions in the limit $t\to\infty$, we insist that $\Re(m^2)\geq{}0$. To assess the boundedness of the solutions~\eqref{eq:usol} in the spatial dependence, we note that the Bessel function $J_n(mr)$ is singular at $r\to\infty$ for all $m\in\mathbb{C}$ such that $|\Im(m)|>0$, and the Hankel function $H^{(2)}_{n}(mr)$ is singular at the origin and at $r\to\infty$ for all $m\in\mathbb{C}$ with $\Im(m)>0$. In addition, in relation to~\eqref{eq:BCs}, we note that $H^{(2)}_{n}(mr)$ vanishes in the limit $r\to\infty$ if~$\Im(m)\leq{}0$. Moreover, $r^n$ (resp. $r^{-n}$) is singular in the limit $r\to\infty$ (resp. $r\to{}0$). On account of their singularity at the origin, $H^{(2)}_{n}$ and $r^{-n}$ are inadmissible in the droplet domain $\Omega_{\DD,0}\supset\{0\}$. Hence, in the droplet domain, it must hold that~$B=0$ and~$F=0$. Conversely, $J_n(mr)$ and $r^n$ are 
inadmissible in the ambient domain, in view of their singularity at $r\to\infty$. Hence, in the ambient domain, it must hold that $A=0$ and $E=0$.

Summarizing, we obtain the following general bounded complex-valued velocity solutions in the droplet and ambient domains:
\begin{subequations}
\label{eq:usol_DA}
\begin{align}
u_{\DD,r}
&=
e^{- \frac{\eta_\DD}{ \rho_\DD} m_\DD^2 t}
\,
\Big[ D_\DD \cos( n_\DD \theta) -  C_\DD \sin( n_\DD \theta) \Big]
\notag\\
&\hspace{100pt}\times
\frac{n_\DD}{r} \Big[ A  J_{n_\DD} ( m_\DD r ) + E r^{n_\DD}  \Big]\,,
\label{eq:usol_Dr}
\\
u_{\DD,\theta}&=
e^{- \frac{\eta_\DD}{ \rho_\DD} m_\DD^2 t}\,
\Big[  C_\DD \cos( n_\DD \theta) + D_\DD \sin( n_\DD \theta) \Big]
\notag\\
&\hspace{100pt}\times
\Big[- A \Big( m_\DD J_{n_\DD-1} ( m_\DD r ) - \frac{n_\DD}{r} J_{n_\DD} ( m_\DD r ) \Big) - E n_\DD r^{n_\DD-1} \Big] \,,
\label{eq:usol_Dtheta}
\\
u_{\AA,r}
&=
e^{- \frac{\eta_\AA}{ \rho_\AA} m_\AA^2 t}\,
\Big[ D_\AA \cos( n_\AA \theta) - C_\AA \sin( n_\AA \theta) \Big]
\notag\\
&\hspace{100pt}\times\frac{n_\AA}{r}\Big[ B  H^{(2)}_{n_\AA} ( m_\AA r ) + F r^{-n_\AA} \Big]\,, 
\label{eq:usol_Ar}
\\
u_{\AA,\theta}&=
e^{- \frac{\eta_\AA}{ \rho_\AA} m_\AA^2 t}\,
\Big[  C_\AA \cos( n_\AA \theta) +  D_\AA \sin( n_\AA \theta) \Big]
\notag\\
&\hspace{100pt}\times
\Big[ - B \Big( m_\AA H^{(2)}_{n_\AA-1} ( m_\AA r ) - \frac{n_\AA}{r} H^{(2)}_{n_\AA} ( m_\AA r ) \Big) + F n_\AA r^{-n_\AA-1}   \Big] \,,
\label{eq:usol_Atheta}
\end{align}
\end{subequations}
subject $\Re(m_i^2)\geq{}0$ $(i\in\{\AA,\DD\})$ and $\Im(m_{\AA})<0$.

%===============================================================================================================%

\subsubsection{Pressure solutions }
To facilitate the derivation of the infinitesimal pressure solutions associated with~\eqref{eq:usol_DA}, we first note that by virtue of~\eqref{eq:sharp_1st_a} and~\eqref{eq:sharp_1st_b}, the pressure solutions are harmonic functions. Considering functions that are sinusoidal and periodic in the angular dependence, that conform to~\eqref{eq:usol_DA} in the temporal dependence, and that are appropriately bounded, we find the following general expression for~$p_{\DD}$:
\begin{equation}
\label{eq:pressure_gen}
p_{\DD}(r,\theta,t)=
e^{- \frac{\eta_\DD}{ \rho_\DD} m_\DD^2 t}\,
r^{\tilde{n}_{\DD}}\big(\tilde{D}_{\DD}\cos(\tilde{n}_{\DD}\theta)+\tilde{C}_{\DD}\sin(\tilde{n}_{\DD}\theta)\big)\,,
\end{equation}
with $\tilde{n}_{\DD}\in\mathbb{Z}_{\geq0}$ and $\tilde{C}_{\DD},\tilde{D}_{\DD}\in\mathbb{C}$ arbitrary constants. By substituting~\eqref{eq:usol_DA} and~\eqref{eq:pressure_gen} 
into~\eqref{eq:sharp_1st_a}, we deduce that compatibility between~\eqref{eq:pressure_gen} and~\eqref{eq:usol_DA} imposes that the index $\tilde{n}_{\DD}\in\mathbb{Z}_{\geq0}$ and coefficients $\tilde{C}_{\DD},\tilde{D}_{\DD}\in\mathbb{C}$ satisfy
\begin{equation}
\tilde{D}_{\DD}=D_{\DD}E\eta_{\DD}m_{\DD}^2,\qquad
\tilde{C}_{\DD}=-C_{\DD}E\eta_{\DD}m_{\DD}^2,\qquad
\tilde{n}_{\DD}=n_{\DD}\,.
\end{equation}
The infinitesimal pressure solution in the ambient domain can be determined similarly. In summary, we obtain:
\begin{subequations}
\label{eq:psol_DA}
\begin{alignat}{2}
p_{\DD} & = 
E \eta_\DD m_\DD^2\, 
e^{- \frac{\eta_\DD}{ \rho_\DD} m_\DD^2 t}\,
r^{n_\DD} \big[ D_\DD \cos(n_\DD \theta) - C_\DD \sin(n_\DD \theta) \big]  \,,\\
p_{\AA} & = - F \eta_\AA m_\AA^2\, 
e^{- \frac{\eta_\AA}{ \rho_\AA} m_\AA^2 t}
r^{-n_\AA} \big[ D_\AA \cos(n_\AA \theta) - C_\AA \sin(n_\AA \theta) \big]  \,.
\end{alignat}
\end{subequations}

\begin{remark}
It is noteworthy that the Bessel and Hankel functions that appear in the velocity solutions~\eqref{eq:usol_DA} are absent in the pressure solutions~\eqref{eq:psol_DA}.
This can be rationalized by noting that these special functions originated directly from the (pressure free) vorticity equation and thus satisfy the momentum equation for a uniform pressure field.
\end{remark}

%===============================================================================================================%

\subsection{Interface conditions}
The general solutions for the pressure and velocity fields in the droplet and ambient domains according to~\eqref{eq:usol_DA} and~\eqref{eq:psol_DA}, involve twelve unknown coefficients: $A$, $B$, $C_\DD$, $C_\AA$, $D_\DD$, $D_\AA$, $E$, $F$, $m_\DD$, $m_\AA$, $n_\DD$ and $n_\AA$. We next extract from the general solutions the subspace that complies with the kinematic interface 
conditions~\eqref{eq:sharp_1st_c} and~\eqref{eq:sharp_1st_d}, and the dynamic condition~\eqref{eq:sharp_1st_e}, by introducing auxiliary conditions on the coefficients.

\subsubsection{Kinematic compatibility conditions}
The kinematic conditions~\eqref{eq:sharp_1st_c}-\eqref{eq:sharp_1st_d} can be equivalently reformulated as:
\begin{subequations} 
\label{eq:kincond}
\begin{alignat}{2}
    &\B{u}_\DD \cdot \B{n}_0 = \B{u}_\AA \cdot \B{n}_0=  \mathcal{V}_1 \qquad&& \textrm{ on } \Gamma_0 \,, \label{eq:kincond_a}\\
    &\B{u}_\DD \cdot \B{t}_0 = \B{u}_\AA \cdot \B{t}_0 && \textrm{ on } \Gamma_0 \,. \label{eq:kincond_b}
\end{alignat}
\end{subequations}
The generating solutions of the interface normal and tangent vectors corresponding to the circular droplet, $\B{n}_0$ and $\B{t}_0$, simply coincide with the radial and angular basis vectors, respectively. The kinematic condition~\eqref{eq:kincond_a} (resp.~\eqref{eq:kincond_b}) thus pertains to the radial (resp. angular) components of~$\B{u}$ in~\eqref{eq:usol_Dr} and~\eqref{eq:usol_Ar} (resp.~\eqref{eq:usol_Dtheta} and~\eqref{eq:usol_Atheta}). To impose~\eqref{eq:kincond_a}, we require the infinitesimal interface velocity $\mathcal{V}_1$ corresponding to~\eqref{eq:R_delta}:
\begin{equation}
\mathcal{V}(\theta,t)
=
\partial_t{}R_{\delta}(\theta,t)    
=
-\delta R_0 \gamma\, e^{- \gamma t}\,\big( \beta \cos ( k \theta ) + \sqrt{1-\beta^2 } \sin ( k \theta ) \big) \,,
\end{equation}
where we retain the entire complex form, to facilitate the exposition. Noting that $\mathcal{V}_0$ vanishes and, hence, $\mathcal{V}=\delta\mathcal{V}_1$, we infer 
from~\eqref{eq:kincond_a} that
\begin{equation}
\label{eq:compatcond}
\begin{aligned}
&e^{- \frac{\eta_\DD}{ \rho_\DD} m_\DD^2 t}\,
\Big[ D_\DD \cos( n_\DD \theta)  -  C_\DD \sin( n_\DD \theta) \Big] \, \Big[ R^{-1}_0 A n_\DD J_{n_\DD} ( m_\DD R_0 ) + E n_\DD R_0^{n_\DD-1}  \Big]  \\
&\qquad= e^{- \frac{\eta_\AA}{ \rho_\AA} m_\AA^2 t}\, \Big[ D_\AA \cos( n_\AA \theta) - C_\AA \sin( n_\AA \theta) \Big]\,\Big[ R^{-1}_0 B n_\AA H^{(2)}_{n_\AA} ( m_\AA R_0 ) + F n_\AA R_0^{-n_\AA-1} \Big]\\
&\qquad\qquad= -R_0 \gamma\, e^{- \gamma t} \, \big( \beta \cos ( k \theta ) + \sqrt{1-\beta^2 }\sin ( k \theta ) \big)\, .
\end{aligned}
\end{equation}
The equalities in~\eqref{eq:compatcond} must hold for all~$\theta\in[0,2\pi)$ and all~$t\in\mathbb{R}_{>0}$. Keeping $\theta$ fixed and varying $t$, one can infer that the temporal exponents must coincide. Subsequently, by fixing~$t$ and varying~$\theta$, it follows that all parameters that characterize the trigonometric terms must be the same. Hence,
\begin{subequations}\label{eq:nDC}
\begin{alignat}{2}
& \frac{\eta_\DD}{\rho_\DD} m_\DD^2 = \frac{\eta_\AA}{\rho_\AA} m_\AA^2 = \gamma \,, \label{eq:mDDmAAgamma}\\ 
& n_\DD = n_\AA = k \,,\\
& D_\DD = D_\AA = \beta \,,\\
& C_\DD = C_\AA = -\sqrt{1-\beta^2} \,.
\end{alignat}
\end{subequations}
In the sequel, we continue to use~$m_\DD$ and~$m_\AA$ in the arguments of the Bessel and Hankel functions, but we tacitly suppose the relation to~$\gamma$ per~\eqref{eq:mDDmAAgamma}. Substitution of~\eqref{eq:nDC} in \eqref{eq:kincond} yields the following three conditions on~$A,B,E$ and~$F$:
\begin{subequations}
\label{eq:kinematic_3}
\begin{align}
A k R^{-1}_0 J_k (m_\DD R_0) + E k R_0^{k-1} &= - \gamma R_0   \,, \label{const1}
\\
B k R_0^{-1} H^{(2)}_{k} ( m_\AA R_0 ) + F k R_0^{-k-1} &= - \gamma R_0  \,, \label{const2}\\
- A \left( m_\DD J_{k-1} ( m_\DD R_0) - k R_0^{-1} J_k ( m_\DD R_0 ) \right) -E k R_0^{k-1}\qquad\qquad& 
\notag\\
+ B \left( m_\AA H_{k-1}^{(2)}(m_\AA R_0) - k R_0^{-1} H_k ^{(2)}(m_\AA R_0 ) \right) - F k R_0^{-k-1} &= 0  \,. \label{const3}
\end{align}
\end{subequations}

%===============================================================================================================%

\subsubsection{Dynamic compatibility conditions}
The dynamic condition~\eqref{eq:sharp_1st_e} can be separated into radial and angular components to form the following two conditions:
\begin{subequations}
\label{eq:dyncond_ab}
\begin{alignat}{2}
- \eta_{\DD}  \B{n}_0 \left(\nabla \B{u}_{\DD} + \left(\nabla \B{u}_{\DD} \right)^T \right) \B{n}_0 + p_{\DD}  + \eta_{\AA}  \B{n}_0 \left(\nabla \B{u}_{\AA} + \left(\nabla \B{u}_{\AA} \right)^T \right) \B{n}_0 - p_{\AA} &= \sigma_{\DD\AA} \kappa_1  
&\quad& \textrm{ on } \Gamma \,, \label{eq:dyncond_a} 
\\
 - \eta_{\DD}  \B{t}_0 \left(\nabla \B{u}_{\DD} + \left(\nabla \B{u}_{\DD} \right)^T \right) \B{n}_0 + \eta_{\AA}  \B{t}_0 \left(\nabla \B{u}_{\AA} + \left(\nabla \B{u}_{\AA} \right)^T \right) \B{n}_0  &= 0  &\quad& \textrm{ on } \Gamma\,. \label{eq:dyncond_b}
\end{alignat}
\end{subequations}
To elaborate on the dynamic interface condition~\eqref{eq:dyncond_a}, we require the first-order perturbation of the interface curvature. From the postulated interface displacement~\eqref{eq:R_delta}, the complex-valued form of the curvature can be derived up to second-order terms:
\begin{multline}
\kappa(\theta,t) 
=
\kappa_0+\delta\kappa_1(\theta,t)+\mathcal{O}(\delta^2)
\\
=
R_0^{-1} + \delta R^{-1}_0 e^{-\gamma t} \big( k^2 - 1 \big) \big( \beta \cos ( k \theta ) - \sqrt{1-\beta^2} \sin ( k \theta) \big)  + \mathcal{O} (\delta^2 ) \,.
\end{multline}
From the expressions for the velocity~\eqref{eq:usol_DA} and pressure~\eqref{eq:psol_DA}, the relation between the coefficients in~\eqref{eq:nDC},
and the dynamic conditions~\eqref{eq:dyncond_ab}, it then follows that
\begin{subequations}
\label{eq:interface_cond_d}
\begin{align}
&- A \eta_\DD \Big[ 2 m_\DD k {R_0}^{-1} J_{k-1}(m_\DD {R_0}) - 2 (k+1) k  {R_0}^{-2} J_k (m_\DD {R_0}) \Big] 
\notag\\
&\qquad+ B \eta_\AA  \Big[ 2 m_\AA k {R_0}^{-1} H_{k-1}^{(2)}(m_\AA {R_0}) - 2 (k+1) k  {R_0}^{-2} H_k ^{(2)}(m_\AA {R_0}) \Big] 
\notag\\
&\qquad+ E\eta_\DD  \Big[ m_\DD^2 {R_0}^k  - 2 (k-1) k {R_0}^{k-2} \Big]+ F \eta_\AA \Big[ m_\AA^2 {R_0}^{-k}  - 2 (k+1) k  {R_0}^{-k-2} \Big] 
\notag\\
&\qquad\qquad= \sigma_{\DD\AA} {R_0}^{-1} (k^2 - 1) \,, \label{const4}
\\
&-A \eta_\DD \Big[ 2 m_\DD {R_0}^{-1} J_{k-1}(m_\DD {R_0}) + \big( m_{\DD}^2 - 2 (k+1) k  {R_0}^{-2} \big) J_k (m_\DD {R_0}) \Big] 
\notag\\
&\qquad+ B \eta_\AA \Big[ 2 m_\AA {R_0}^{-1} H_{k-1}^{(2)}(m_\AA {R_0}) + \left( m_\AA^2 - 2 (k+1) k  {R_0}^{-2} \right) H_k ^{(2)}(m_\AA {R_0}) \Big] 
\notag\\
&\qquad+ 2 E \eta_\DD (k-1) k {R_0}^{k-2}  + 2 F \eta_\AA (k+1) k  {R_0}^{-k-2}  
\notag\\
&\qquad\qquad= 0 \,. \label{const5}
\end{align}
\end{subequations}

%===============================================================================================================%

\subsection{Dispersion relation}

For each wave number~$k\in\mathbb{N}_{\geq{}2}$, the corresponding characteristic temporal response coefficient of the solution, $\gamma\in\mathbb{C}$, as well as the mode shapes, encoded in the remaining free parameters, follow from a solution-existence condition. To elucidate this condition, we first recall that the general bounded complex-valued velocity and pressure solutions of the partial-differential equations \eqref{eq:sharp_1st_a}-\eqref{eq:sharp_1st_b}, subject to the limit condition~\eqref{eq:BCs}, are given by~\eqref{eq:usol_DA} and~\eqref{eq:psol_DA}. These general solutions contain twelve coefficients. Eight of these coefficients are determined by the kinematic interface condition~\eqref{eq:sharp_1st_c}, in accordance with~\eqref{eq:nDC}. The kinematic conditions~\eqref{eq:sharp_1st_c} and~\eqref{eq:sharp_1st_d} imply that the remaining four coefficients, $A,B,E$ and~$F$ must satisfy the three identities in~\eqref{eq:kinematic_3}. The dynamic condition~\eqref{eq:sharp_1st_e} demands that, in addition, these four coefficients satisfy the two identities in~\eqref{eq:interface_cond_d}. The remaining five conditions on the coefficients can be cast in the form
\begin{equation}
\label{eq:sysofeq}
\underbrace{%
\begin{pmatrix}
{a}_{11}&0&{a}_{13}&0
\\
0&{a}_{22}&0&{a}_{24}
\\
{a}_{31}&{a}_{32}&{a}_{33}&{a}_{34}
\\
{a}_{41}&{a}_{42}&{a}_{43}&{a}_{44}
\\
{a}_{51}&{a}_{52}&{a}_{53}&{a}_{54}
\end{pmatrix}}_{\B{A}(k,\gamma)}
\begin{pmatrix}
A \\
B \\
E \\
F 
\end{pmatrix} = 
\underbrace{%
\begin{pmatrix}
-\gamma R_0 \\
-\gamma R_0 \\
0 \\
\sigma_{\DD\AA} R^{-1}_0 (k^2 - 1) \\
0
\end{pmatrix}}_{{\B{b}}(k,\gamma)}  \,,
\end{equation}
in such a manner that the first three equations in~\eqref{eq:sysofeq} represent~\eqref{eq:kinematic_3} and the latter two represent~\eqref{eq:interface_cond_d}. Noting the dependence of~\eqref{eq:kinematic_3} and~\eqref{eq:interface_cond_d} on the wave number~$k$, and recalling the dependence 
of~$m_{\DD}$ and~$m_{\AA}$ in these equation on the temporal response coefficient~$\gamma$ via~\eqref{eq:mDDmAAgamma}, we infer that the entries of~${\B{A}}$ depend on~$k$ and~$\gamma$. With five constraints and four unknowns,
the system of equations~\eqref{eq:sysofeq} is formally over-constrained, and a solution is non-existent unless the right-hand-side vector~${\B{b}}(k,\gamma)$ is in the column space of~${\B{A}}(k,\gamma)$. The relation between the existence of a solution and the condition
\begin{equation}
\label{eq:exist_cond}
{\B{b}}(k,\gamma)\in\operatorname{span}(\operatorname{col}({\B{A}}(k,\gamma)))     \,,
\end{equation}
is indicative of the fact that only specific combinations of the wave number $k\in\mathbb{N}_{\geq{}2}$ and the temporal response coefficient $\gamma\in\mathbb{C}$ in the postulated interface configuration~\eqref{eq:Gamma_delta} correspond to a natural response of the droplet. 

To determine the combinations $(k,\gamma)$ for which the existence condition~\eqref{eq:exist_cond} is fulfilled, we note that~\eqref{eq:exist_cond} is equivalent to
\begin{equation}
\label{eq:det_cond}
\det\big((\B{A}\mid\B{b})(k,\gamma)\big)=0    \,,
\end{equation}
where $(\B{A}\mid\B{b})$ corresponds to $\B{A}$ augmented by~$\B{b}$. The equivalence between~\eqref{eq:exist_cond} and~\eqref{eq:det_cond} follows from the fact that the column vectors of~$\B{A}$ are linearly independent for all $(k,\gamma)$ and, hence, the augmented matrix is singular if and only if the vector~$\B{b}$ resides in the column space of~$\B{A}$. Moreover, by virtue of the linear independence of the columns of~$\B{A}$, if~\eqref{eq:det_cond} holds, then~\eqref{eq:sysofeq} has a unique solution. This solution corresponds to the coefficients~$(A,B,E,F)$ that, in combination with~\eqref{eq:nDC}, define the droplet and ambient velocity-pressure pairs corresponding to $(k,\gamma)$ according to~\eqref{eq:usol_DA} and~\eqref{eq:psol_DA}.

To facilitate and generalize the root-finding of the determinant in~\eqref{eq:det_cond}, we non-dimension\-alize the matrix entries of the augmented matrix based on the droplet density, $\rho_{\DD}$, droplet viscosity, $\eta_{\DD}$ and droplet radius~$R_0$. The non-dimensionalized parameters are indicated with a tilde diacritic. Additionally, we introduce the following condensed notation:
\begin{alignat*}{6}
&\mathcal{J} = J_k ( \tilde{m}_\DD (\tilde{\gamma}) )\,,\qquad && \mathcal{H} = H^{(2)}_k ( \tilde{m}_\AA (\tilde{\gamma}) ) \,,\qquad && \zeta = 2(k-1)k \,, \\
&\hat{\mathcal{J}} =  \tilde{m}_\DD (\tilde{\gamma}) J_{k-1} ( \tilde{m}_\DD (\tilde{\gamma}) )  \,,\qquad && \hat{\mathcal{H}} = \tilde{m}_\AA (\tilde{\gamma}) H^{(2)}_{k-1} ( \tilde{m}_\AA (\tilde{\gamma}) ) \,,\qquad && \xi = 2 (k+1) k \,.
\end{alignat*}
The non-dimensionalized augmented matrix can then be expressed as
\begin{multline}
\label{eq:Abtilde}
(\tilde{\B{A}}|\tilde{\B{b}})(k,\tilde{\gamma})=
\\
\begin{pmatrix}
k \mathcal{J} & 0 & k & 0 & - \tilde{\gamma} \\
0 & k\mathcal{H} & 0 & k & - \tilde{\gamma}  \\
- \hat{\mathcal{J}} + k \mathcal{J} & \hat{\mathcal{H}} - k \mathcal{H} & -k & -k & 0 \\
-2 k \hat{\mathcal{J}} + \xi \mathcal{J} & \tilde{\eta}_\AA \big[ 2 k \hat{\mathcal{H}} - \xi \mathcal{H} \big] & \tilde{m}^2_\DD - \zeta & \tilde{\eta}_\AA \big [ \tilde{m}^2_\AA - \xi \big] & \tilde{\sigma}_{\DD\AA}\left(k^2 - 1\right)  \\
-2\hat{\mathcal{J}} - \left(\tilde{m}^2_\DD - \xi \right) \mathcal{J} & \tilde{\eta}_\AA \big[ 2 \hat{\mathcal{H}} + \left(\tilde{m}^2_\AA - \xi \right) \mathcal{H}\big] & \zeta & -\tilde{\eta}_\AA \xi & 0
\end{pmatrix}.
\end{multline}
It is not generally feasible to determine the roots of~$\det((\tilde{\B{A}}|\tilde{\B{b}})(k,\tilde{\gamma}))$ with respect to~$\tilde{\gamma}$ in closed form and, in practice, it is necessary to revert to a numerical root-finding algorithm. Once a root has been determined, one can extract the kernel of the augmented matrix~\eqref{eq:Abtilde} and scale the corresponding vector such that its fifth entry is minus one, to obtain the coefficients $\tilde{A},\tilde{B},\tilde{E},\tilde{F}$.

\begin{remark}
The roots of~$\det((\tilde{\B{A}}|\tilde{\B{b}})(k,\tilde{\gamma}))$ are not unique: one
can infer that
\begin{equation}
\big[(\tilde{\B{A}}|\tilde{\B{b}})(k,\tilde{\gamma}{}^*)\big]
=
\big[(\tilde{\B{A}}|\tilde{\B{b}})(k,\tilde{\gamma})\big]^*\,,
\end{equation}
where $(\cdot)^*$ denotes complex conjugation. Because the eigenvalues of the complex conjugate of a matrix are the complex conjugates of the original eigenvalues, it follows that if $\tilde{\gamma}$ is a root of~$\det((\tilde{\B{A}}|\tilde{\B{b}})(k,\tilde{\gamma}))$, then so is~$\tilde{\gamma}{}^*$. Since, in addition, it must hold that $\Re(\tilde{\gamma})>0$, it suffices to consider roots in the fourth quadrant of the complex plane. Noting that the entries of the augmented matrix are analytic functions, one can infer that so is its determinant. This implies that the roots of~$\det((\tilde{\B{A}}|\tilde{\B{b}})(k,\tilde{\gamma}))$ form a totally disconnected set and, accordingly, for each root there exists a neighborhood in which that root is unique. A detailed investigation of the uniqueness of the roots of~$\det((\tilde{\B{A}}|\tilde{\B{b}})(k,\tilde{\gamma}))$ in the fourth quadrant is beyond the scope of this work. In our numerical root-finding procedure, we have verified that there are no other roots in a region around the found root.
\end{remark}

By virtue of the complex representation of the interface parametrization~\eqref{eq:R_delta}, we obtain the real-valued velocity and pressure fields by taking the real parts of~\eqref{eq:usol_DA} and \eqref{eq:psol_DA} after substitution of the relations \EQ{nDC}, and the temporal response coefficient~$\gamma$ and the corresponding coefficients $A, B, E$, and $F$. Table~\ref{tab:solutionfieldparameters} provides computed parameter values for the physical setting outlined in Table~\ref{tab:model_param_vals}, representing a water-in-air picoliter-sized droplet. For completeness, we mention that we have applied
\emph{Mathematica}'s root-finder to determine $\tilde{\gamma}:=\tilde{\gamma}_k$ in the fourth quadrant of the complex plane such that $\det ((\tilde{\B{A}}| \tilde{\B{b}} )(k,\tilde{\gamma}_k))=0$. The dimensions of parameters~$E$ and~$F$ depend on the mode number~$k$. As a result, as~$k$ increases, the values of~$E$ and~$F$ grow rapidly, conveying that these parameters are ill-conditioned in terms of~$k$. Furthermore, as the mode number~$k$ increases, the frequency and damping rate of the corresponding oscillation, both encoded in~$\gamma$, increase. This implies that if a droplet sustains an initial perturbation that is characterized by multiple modes, the higher wave-number modes decay quickly, and low order modes dominate the long-term dynamics of viscous-in-viscous oscillating droplets.

\begin{table}[]\centering
\caption{Modal solution parameter values for a water droplet of radius $R= \sqrt{2}\!\times\!10^{1}\,\textrm{\textmu m}$ suspended in air.} \label{tab:solutionfieldparameters}
\footnotesize
\begin{tabular}{l|ccccc}
$k$ & $\gamma\,\, [ \textrm{ s}^{-1} ]$     & $A\,\, [10^{-11} \textrm{ m}^2\textrm{s}^{-1}]$ & $B\,\, [10^{-4}\textrm{ m}^2\textrm{s}^{-1}]$ & $E\,\, [10^{5k-7}\textrm{ m}^{2-k}\textrm{s}^{-1}]$ & $F\,\, [10^{-5 k - 4}\textrm{ m}^{2+k}\textrm{s}^{-1}]$ \\[0.2cm]\hline\hline\\[-0.2cm]
2   & $18788.18393 $    & $4664.160935$     & $1.142101737$        & $2.618366811$     & $1.904600443$    \\
    & $-390396.1271\,i$ & $+137.5420287\,i$ & $+4.485396518\,i$    & $+194.0864949\,i$ & $+1.883632636\,i$\\[0.2cm]
3   & $53722.95262$     & $376.1242266$     & $7.390256000$        & $3.392996002$     & $4.168582183$    \\
    & $-777097.6733\,i$ & $-953.0487542\,i$ & $+6.296111161\,i$    & $+182.1390754\,i$ & $+3.618449771\,i$\\[0.2cm]
4   & $104333.5005$     & $-123.9009766$    & $14.90136892$        & $3.548390076$     & $7.942104534$    \\
    & $-1223261.977\,i$ & $-195.0741019\,i$ & $+ 1.040690443\,i$   & $+152.0972759\,i$ & $+6.209497220\,i$\\[0.2cm]
5   & $170139.1932$     & $-53.99114858$    & $16.35189970$        & $3.366783037$     & $14.10551086$    \\
    & $-1722774.521\,i$ & $-11.92064586\,i$ & $-10.42851063\,i$    & $+121.2100680\,i$ & $+10.15645338\,i$\\[0.2cm]
6   & $250754.7829$     & $-12.54654878$    & $8.074025870$        & $3.017014114$     & $23.99270873$    \\
    & $-2270070.007\,i$ & $+6.299980519\,i$ & $-20.94948730\,i$    & $+94.14089637\,i$ & $+16.16915562\,i$\\[-0.2cm]
\end{tabular}
\end{table}

%===============================================================================================================%
%===============================================================================================================%
%===============================================================================================================%

\section{Numerical experiments}
\label{sec:num_exp}
The free-boundary problem~\eqref{eq:sharp_strong} formally represents the 
sharp-interface limit of the Abels--Garcke--Gr\"un NSCH model~\eqref{eq:strong}, provided that the mobility is appropriately scaled in the limit $\varepsilon\to{}+0$.
For sufficiently small~$\delta$, the oscillating-droplet solutions derived in Section~\ref{sec:deriv} can therefore serve to investigate the approach of the diffuse-interface solution to the sharp-interface limit solution. In this section, we
investigate this sharp-interface limit numerically, by means of an adaptive finite-element method. Specifically, we focus on the scaling of the mobility parameter $m:=m_{\varepsilon}$ in the limit $\varepsilon\to{}+0$, and investigate the deviation of the diffuse-interface solution from the sharp-interface solution in relation to~$m$.
As reference solutions, we consider the lowest mode of oscillation ($k=2$), as well as the next higher doubly symmetric mode ($k=4$), for a viscous droplet in a viscous ambient with parameter values according to Table~\ref{tab:model_param_vals}.
The corresponding coefficients of the velocity solution~\eqref{eq:usol_DA} and 
pressure solution~\eqref{eq:psol_DA} fields are presented 
in Table~\ref{tab:solutionfieldparameters}.

\begin{table}[!b]\centering
\caption{Physical and numerical parameter values of the considered numerical experiments. Entries marked with the symbol~$\ast$ indicate a range of values, which will be specified in the text.}
\label{tab:model_param_vals}
\begin{tabular}{cc|cc|ccc|cccc}
\multicolumn{2}{c|}{Droplet}&\multicolumn{2}{c|}{Ambient}&\multicolumn{3}{c|}{Interface}&\multicolumn{4}{c}{Numerical approximation}\\
$\rho_{\DD}$ & $\eta_{\DD}$ & $\rho_{\AA}$  & $\eta_{\AA}$  & $\sigma_{\DA}$  & $\varepsilon$       & $m$       & $\tau$             & $h_0$ & $L_{\MAX}$ & $K$ \\[0.1cm]
\hline
\hline & & & & & & &  & &  & \\[-0.3cm]
$\frac{\textrm{kg}}{\textrm{m}^d}$ &$\frac{\textrm{kg}\,\textrm{m}^{2-d}}{\textrm{s}}$& $\frac{\textrm{kg}}{\textrm{m}^d}$  & $\frac{\textrm{kg}\,\textrm{m}^{2-d}}{\textrm{s}}$ & $\frac{\textrm{kg}\,\textrm{m}^{3-d}}{\textrm{s}^2}$ & m & $\frac{\textrm{m}^d\,\textrm{s}}{\textrm{kg}}$ & $10$ \textmu s & \textmu m & --- & --- \\[0.2cm]
$10^3$ & $10^{-3}$ & $1$ & $1.813\!\times\!10^{-5}$ & $7.28\!\times\!10^{-2}$ & $\ast$ & $\ast$ & $2^{-7}$ & $5$ & $\ast$ & $\ast$ \\[0.1cm]
\end{tabular}
\end{table}

\subsection{Setup and discretization}
The oscillating-droplet test cases that we consider pertain to doubly symmetric modes. The setup of the test cases is similar to that in~\cite[Sec.5]{Demont:2022dk}. To reduce computational expense, we exploit the symmetry of the configurations and consider only one quarter of the droplet-ambient domain.
We regard a domain~$\Omega=(0,50)^2\,\mu\text{m}^2$ and prescribe symmetry conditions on $\Gamma_{\text{sym}}:=\{(x_1,x_2)\in\partial\Omega:\{x_1=0\}\cup\{x_2=0\}\}$. Because the linear sharp-interface solution is in fact defined on the generating circular droplet domain and the corresponding ambient domain according to~\eqref{eq:gensol_a}, while the diffuse-interface model exhibits a moving interface, we prescribe auxiliary conditions in accordance with an initially circular droplet. Specifically, with reference to~\eqref{eq:R_delta}, we select~$t_0$ such that $-\Im(\gamma)t_0=\pi/2$ and, hence, $R_{\delta}(\theta,t_0)=R_0$, and prescribe initial data corresponding to the reference solution at~$t_0$ and boundary data corresponding to $t+t_0$. The complementary part of the boundary, $\Gamma_{\text{ext}}:=\partial\Omega\setminus\Gamma_{\text{sym}}$, is furnished with Dirichlet conditions for velocity and homogeneous Neumann conditions for the order parameter and the chemical potential:
\begin{equation}
\left.
\begin{aligned}
\B{u}(\cdot,t)&=\delta\B{u}_{\AA}(\cdot,t_0+t)    
\\
\partial_n\varphi&=0
\\
\partial_n\mu&=0
\end{aligned}    
\right\}\quad\text{on }\Gamma_{\text{ext}}, \text{ for }t\in[0,T)\,,
\end{equation}
where $\delta\B{u}_{\AA}$ corresponds to the ambient velocity solution~\eqref{eq:usol_DA} with appropriate coefficients and scaling~$\delta$,
and $[0,T)$ denotes the time-interval under consideration. For the aforementioned combination of boundary conditions, the pressure variable $p$ is only determined up to a constant. We impose the auxiliary condition that~$p$ vanishes on average. 

We impose an initial condition for the order parameter corresponding to a circular interface, in accordance with the initial configuration of the sharp-interface reference solution, viz.
\begin{equation}
\label{eq:phi0}
\varphi(\B{x},0)=\varphi_0(\B{x}):=\tanh\bigg(\frac{d_{\pm}(\B{x},\Gamma_0)}{\sqrt{2}\varepsilon}\bigg)\,,    
\end{equation}
where $d_{\pm}(\B{x},\Gamma_0)$ represents the signed distance from~$\B{x}$ to~$\Gamma_0$. The function $s\mapsto\tanh(s/\sqrt{2}\varepsilon)$ corresponds to an equilibrium solution of the Cahn--Hilliard equations for the phase field 
in one spatial dimension and, accordingly, the phase field~\eqref{eq:phi0} is meta-stable if~$\varepsilon$ is sufficiently small compared to the radius of curvature of~$\Gamma_0$. In conjunction with~\eqref{eq:phi0}, we impose the following initial condition for velocity:
\begin{equation}
\label{eq:u0}
\B{u}(\B{x},0)=
\begin{cases}
\delta\B{u}_{\DD}(\B{x},t_0)&\quad\text{if }\B{x}\in\Omega_{\DD,0}\,,
\\
\delta\B{u}_{\AA}(\B{x},t_0)&\quad\text{if }\B{x}\in\Omega_{\AA,0}\,,
\end{cases}
\end{equation}
where the data in the right member of~\eqref{eq:u0} corresponds to the 
the velocity solutions according to~\eqref{eq:usol_DA} in the droplet and ambient domains. We select the perturbation magnitude~$\delta=10^{-2}$, after verifying that this choice renders the linearization error negligible in comparison to the deviation between the diffuse-interface and the sharp-interface solutions, for the~$\varepsilon$ considered below. Hence, the selected value of~$\delta$ is suitable for our investigation of the sharp-interface limit. The characteristic parameters pertaining to the 
droplet and ambient fluids, and to the interface are reported in Table~\ref{tab:model_param_vals}.

To perform the numerical simulations, we make use of the adaptive finite-element approximation method presented in~\cite{Demont:2022dk}. For coherence, we present a concise overview of the numerical methodology. The weak form of the NSCH equations~\eqref{eq:strong} is discretized with respect to the spatial dependence with $P3-P2$ (Taylor-Hood) $\mathcal{C}^0$ truncated hierarchical B-splines (see \cite{Hughes:2005it, Cottrell:2009ad, Giannelli:2012rr}) for the velocity and pressure fields, and $P3$ $\mathcal{C}^0$ truncated hierarchical B-splines for the order parameter and chemical potential; see~\cite[\S3.1]{Brummelen:2021aw} for further details. The adaptive-refinement procedure is guided by a two-level hierarchical a-posteriori error estimate, and follows the standard SEMR (\texttt{Solve} $\rightarrow$ \texttt{Estimate} $\rightarrow$ \texttt{Mark} $\rightarrow$ \texttt{Refine}) process~\cite{Bertoluzza:2012kx, Dorfler:1996uq}. To improve the robustness of the solution procedure on the coarse meshes that occur in the sequence of adaptive refinements within each time step, an $\varepsilon$-continuation process is introduced, in which the thickness parameter~$\varepsilon$ (and, in conjunction, the mobility~$m \propto \varepsilon^3$) is enlarged for the first~$K$ iterations of the adaptive refinement process; see~\cite{Demont:2022dk, Brummelen:2021aw} for details.
In each time step, the fluid domain is initially covered with a uniform mesh comprising $10\times10$ elements, corresponding to an initial mesh width~$h_0=5\,\mu\text{m}$, and we perform $L_\MAX$ refinement steps. Refinement steps $L=0, 1, \ldots, K-1$ make use of the $\varepsilon$- and $m$-continuation process, while in refinement steps $L=K, \ldots, L_\MAX$ the original parameter values for $\varepsilon$ and $m$ are used. A skew-symmetric formulation according to~\cite{Layton:2008fk} is used for the convective term in the Navier--Stokes equations, enhancing the stability of the discrete approximation by eliminating potential artificial energy production due to deviations from solenoidality in pure-species regions. On the coarse meshes, a first order Backward Euler scheme with second order contractive-expansive splitting of the double-well potential with stabilization~\cite{Wu:2014tg} is employed. On the finest mesh, a second order Crank--Nicolson scheme is applied with implicit treatment of the double-well potential. The second order Crank--Nicolson scheme provides significant better accuracy than the Backward Euler scheme; cf.\ e.g.~\cite{John:2006qi}. For the temporal discretization, we employ a time-step size $\tau = 2^{-7} \times 10\,\mu$s. The parameter setting of the numerical procedure is also summarized in Table~\ref{tab:model_param_vals}. The nonlinear algebraic systems corresponding to the discretized NSCH equations, are solved with a Newton procedure, in which the linear tangent problems are solved with GMRES with a preconditioner based on a partition of the NSCH system into NS and CH subsystems; see~\cite{Demont:2022dk} for further details.

\begin{figure}
\centering
\includegraphics[width=\textwidth]{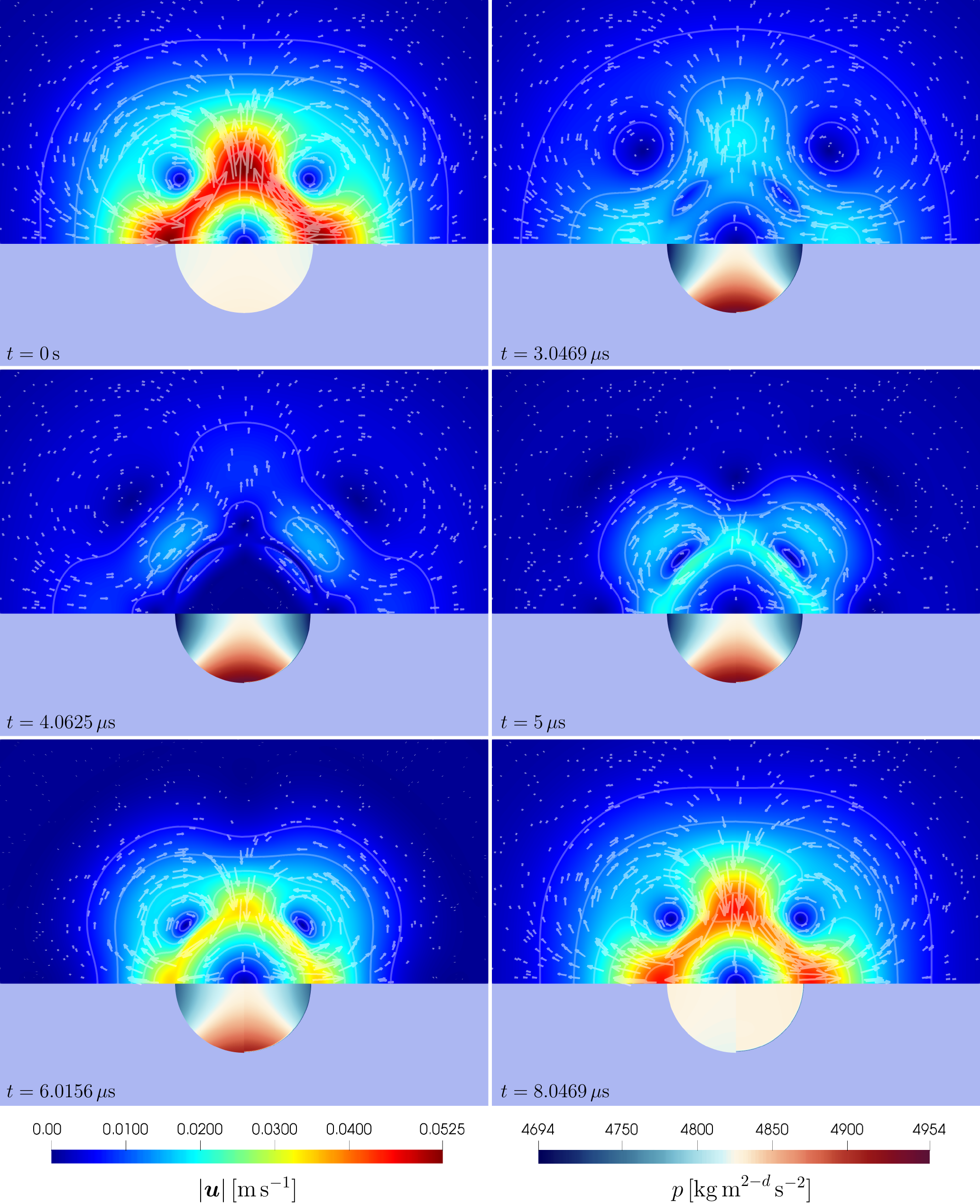}
\caption{Snapshots of the magnitude of the velocity field $|\B{u}|$ (top) and pressure field $p$ (bottom) throughout half a droplet oscillation of mode $k=2$. The left half of each panel displays the analytical sharp-interface solution, while the right half displays the numerical diffuse-interface solution.}
\label{fig:solution_snapshots_k2}
\end{figure}

\begin{figure}
\centering
\includegraphics[width=\textwidth]{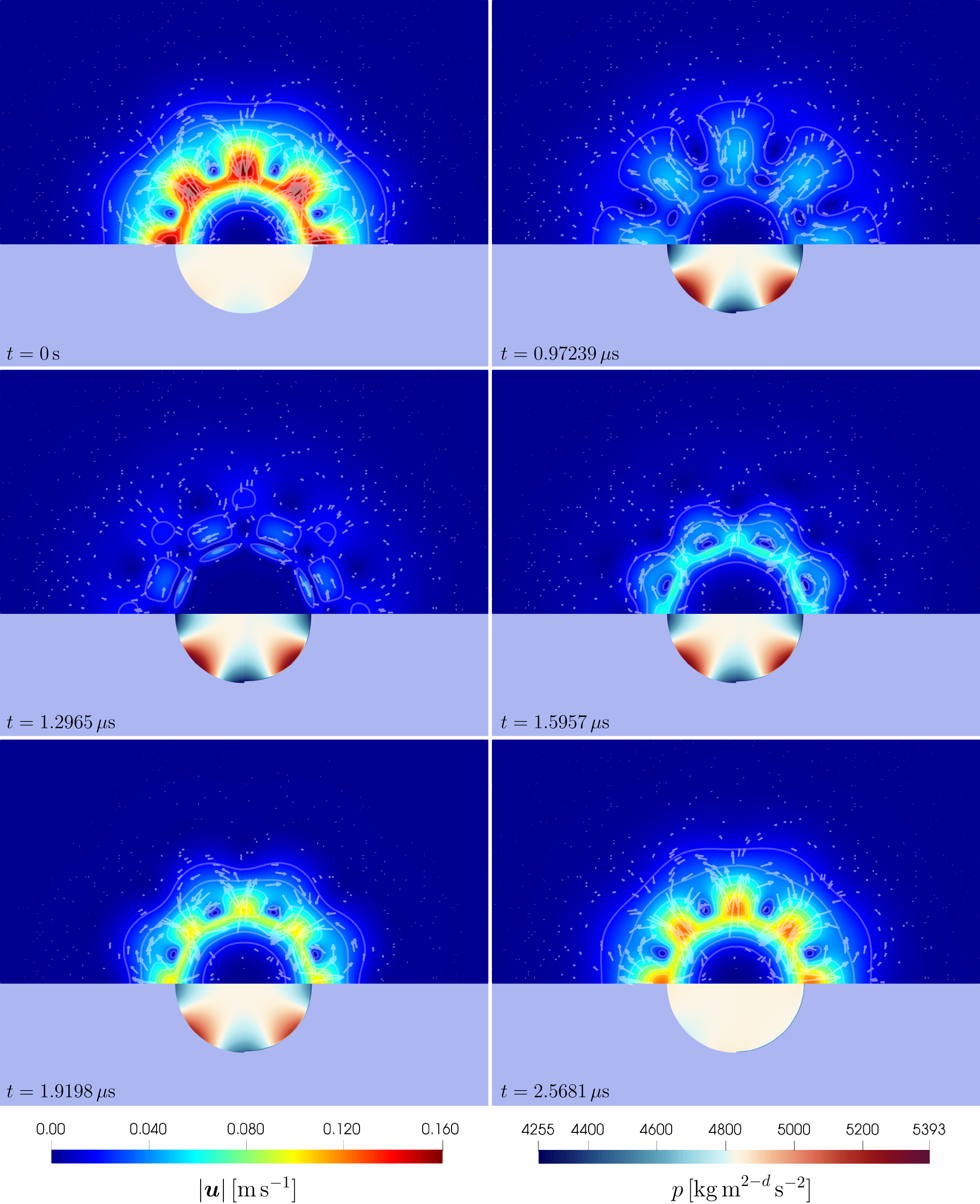}
\caption{Snapshots of the magnitude of the velocity field $|\B{u}|$ (top) and pressure field $p$ (bottom) throughout half a droplet oscillation of mode $k=4$. The left half of each panel displays the analytical sharp-interface solution, while the right half displays the numerical diffuse-interface solution.}
\label{fig:solution_snapshots_k4}
\end{figure}

To illustrate the setup of the numerical experiments, and the resemblance between the analytic sharp-interface solution and the numerical approximation of the diffuse-interface solution for sufficiently small~$\varepsilon$, we conduct numerical experiments with interface-thickness parameter~$\varepsilon=2^{-10}\times{}10^2\,\text{\textmu m}$, mobility~$m=9.5272\times{}10^{-13}\,\text{m}^2\text{s}/\text{kg}$, maximum number of refinement levels~$L_{\MAX}=7$, and number of continuation 
levels~$K=5$. Figures~\FIG{solution_snapshots_k2} and~\FIG{solution_snapshots_k4} display snapshots of the velocity field and pressure field at six time instants. The top (resp.\ bottom) half of each panel displays the velocity (resp.\ pressure) field. The right (resp.\ left) half of each panel depicts the diffuse-interface simulation (resp.\ sharp-interface solution). The figures convey that the sharp-interface solutions and the diffuse-interface solutions are visually indistinguishable.

\subsection{Optimal mobility scaling}
\label{ssec:numexpmobility}
To elucidate the dependence of the diffuse-interface solution in the sharp-interface limit $\varepsilon\to+0$ on the scaling of the mobility $m:=m_{\varepsilon}$, we conduct numerical experiments for a range of combinations of~$\varepsilon$ and~$m$. For each combination of~$\varepsilon$ and~$m$, we determine the deviation relative to the sharp-interface solution according to:
\begin{equation}
\label{eq:devem}
\operatorname{dev}(\varepsilon,m) = \frac{\varnorm{\B{u}^{\varepsilon,m}-\B{u}}_{\Omega\times(0,T)}}{\varnorm{\B{u}}_{\Omega\times(0,T)}} \quad \text{with} \quad \varnorm{\B{u}}_{\Omega\times(0,T)} = \frac{1}{T} \int\limits_{0}^T \norm{\B{u}(\cdot,t)}_{L^2(\Omega)} \,\text{d}t \,,
\end{equation}
where the considered length of the time interval,~$T$, corresponds to half a period of oscillation. We regard a set of decreasing interface thickness parameters $\varepsilon\in\mathscr{E}:=\{2^0,\ldots,2^{-3}\}\varepsilon_\text{max}$ relative to the baseline interface thickness $\varepsilon_\text{max} = 2^{-7}\!\times\!10^{2}\,\textrm{\textmu m} = 7.8125\!\times\!10^{-1}\,\textrm{\textmu m}$. The baseline interface thickness parameter corresponds to approximately 5\% of the droplet radius. For each~$\varepsilon$, we consider mobility parameters in (a relevant subset of) the set $m\in\mathscr{M}:=\{2^{0},2^{-1},\ldots,2^{-12}\}m_\text{max}$
with $m_{\text{max}} = 2.4389632\times10^{-10}\,\textrm{m}^d\textrm{s}/\textrm{kg}$. The range of mobility parameters has been determined empirically such that $[2^{-12},1]m_{\text{max}}$ includes the optimal mobility, i.e.\ the one for which $\operatorname{dev}(\varepsilon,m)$ is minimal, for all $\varepsilon\in\mathscr{E}$. It is to be noted that $\mathcal{E}\times\mathcal{M}$ contains various monomial scalings of the mobility with respect to the interface thickness, viz.\ $m\propto\varepsilon^l$ with $l\in\{0,1,2,3\}$. 

Tables~\ref{tab:mobsk2} and~\ref{tab:mobsk4} present the deviations $\operatorname{dev}(\varepsilon,m)$ for the two modes of oscillation, 
$k=2$ and~$k=4$, respectively. For each~$\varepsilon\in\mathcal{E}$, the entry corresponding to the mobility $m\in\mathcal{M}$ that yields the smallest deviation, is highlighted. One can observe that, indeed, the mobility corresponding to the minimal deviation decreases as~$\varepsilon$ decreases. More precisely, for both modes, the optimal scaling of the mobility parameter with the interface-thickness parameter appears to lie between~$m\propto{}\varepsilon$ and~$m\propto{}\varepsilon^2$. One may moreover note that the entries corresponding to the optimal mobility decrease by a factor of approximately two if~$\varepsilon$ is halved, which indicates that for the considered droplet-oscillation case, the diffuse-interface solution approaches the sharp-interface solution at rate~$\mathcal{O}(\varepsilon)$, provided that the mobility in the diffuse-interface model is appropriately scaled.

\begin{table}[!t]
\centering
\caption{Deviation between the diffuse-interface solution and the sharp-interface solution according to~\eqref{eq:devem}, for $k=2$ for one half period of oscillation.}
\label{tab:mobsk2}
\begin{tabular}{|c|cccc|}
\cline{1-5}
\multicolumn{1}{|c|}{$m/m_\MAX$ } & \rule{0pt}{12pt}$2^{-3}$ \normalsize$\varepsilon_\MAX$ & $2^{-2}$ \normalsize$\varepsilon_\MAX$ & $2^{-1} $ \normalsize$\varepsilon_\MAX$ & $\varepsilon_\MAX$ \\
\cline{1-5}
\rule{0pt}{12pt}$2^{-12}$ & $6.8295\!\times\!10^{-2}$ &  &                           &                          \\
$2^{-11}$ &  &                           &                           &                          \\
$2^{-10}$ &  &                           &                           &                          \\
$2^{-9\phantom{0}}$ & $1.2634\!\times\!10^{-2}$ & $7.3709\!\times\!10^{-2}$ &                           &                          \\
$2^{-8\phantom{0}}$ & \cellcolor{green!20}$8.3534\!\times\!10^{-3}$& $4.4566\!\times\!10^{-2}$ &        &                          \\
$2^{-7\phantom{0}}$ & $1.2941\!\times\!10^{-2}$ & $2.1360\!\times\!10^{-2}$ & $1.0190\!\times\!10^{-1}$ & $1.3056\!\times\!10^{-1}$\\
$2^{-6\phantom{0}}$ & $2.3777\!\times\!10^{-2}$ & \cellcolor{green!20}$1.8713\!\times\!10^{-2}$ & $6.6942\!\times\!10^{-2}$ & $1.2485\!\times\!10^{-1}$\\
$2^{-5\phantom{0}}$ &                           & $4.2245\!\times\!10^{-2}$ & \cellcolor{green!20}$3.5126\!\times\!10^{-2}$ & $1.0979\!\times\!10^{-1}$\\
$2^{-4\phantom{0}}$ &                           &                           & $6.2628\!\times\!10^{-2}$ & $7.9845\!\times\!10^{-2}$\\
$2^{-3\phantom{0}}$ & $1.5565\!\times\!10^{-1}$ & $1.6940\!\times\!10^{-1}$ & $1.5228\!\times\!10^{-1}$ & \cellcolor{green!20}$6.1487\!\times\!10^{-2}$\\
$2^{-2\phantom{0}}$ &                           &                           &                           & $1.7640\!\times\!10^{-1}$\\
$2^{-1\phantom{0}}$ &                           &                           &                           & $4.2783\!\times\!10^{-1}$\\
$1\phantom{^{0-1}}$ &                           &                           &                           & $8.2210\!\times\!10^{-1}$\\
\cline{1-5} 
\end{tabular}
\end{table}

\begin{table}[!t]\centering
\caption{Deviation between the diffuse-interface solution and the sharp-interface solution according to~\eqref{eq:devem}, for $k=4$ for one half period of oscillation.}
\label{tab:mobsk4}
\begin{tabular}{|c|cccc|}
\cline{1-5}
\multicolumn{1}{|c|}{$m/m_\MAX$ } & \rule{0pt}{12pt}$2^{-3}$ \normalsize$\varepsilon_\text{max}$ & $2^{-2}$ \normalsize$\varepsilon_\text{max}$ & $2^{-1}$ \normalsize$\varepsilon_\text{max}$ & $\varepsilon_\text{max}$ \\
\cline{1-5}
\rule{0pt}{12pt}$2^{-11}$ & $1.2922\!\times\!10^{-1}$ & $ $ & $ $ & \\
$2^{-10}$ & $ $ & $ $ & $ $ & \\
$2^{-9\phantom{0}}$ & $ $ & $ $ & $ $ & \\
$2^{-8\phantom{0}}$ & $2.7543\!\times\!10^{-2}$ & $1.4419\!\times\!10^{-1}$ &    & \\
$2^{-7\phantom{0}}$ & \cellcolor{green!20}$2.3628\!\times\!10^{-2}$ & $9.3099\!\times\!10^{-2}$ & $ $ & $ $\\
$2^{-6\phantom{0}}$  & $3.2089\!\times\!10^{-2}$ & $5.5083\!\times\!10^{-2}$ & $ $ & $ $\\
$2^{-5\phantom{0}}$  & $5.2591\!\times\!10^{-2}$ & \cellcolor{green!20}$5.2646\!\times\!10^{-2}$ & $1.5428\!\times\!10^{-1}$ & $ $\\
$2^{-4\phantom{0}}$  &    & $9.4735\!\times\!10^{-2}$ & \cellcolor{green!20}$1.0710\!\times\!10^{-1}$ & $2.9576\!\times\!10^{-1}$\\
$2^{-3\phantom{0}}$   &    &    & $1.2054\!\times\!10^{-1}$ & $2.7358\!\times\!10^{-1}$\\
$2^{-2\phantom{0}}$   & $3.3162\!\times\!10^{-1}$ & $3.5898\!\times\!10^{-1}$ & $2.7331\!\times\!10^{-1}$ & \cellcolor{green!20}$2.0927\!\times\!10^{-1}$\\
$2^{-1\phantom{0}}$   &    &    &    & $2.5501\!\times\!10^{-1}$\\
$1\phantom{^{0-1}}$          &    &    & $ $ & $6.2529\!\times\!10^{-1}$\\
\cline{1-5} 
\end{tabular}
\end{table}

\begin{figure}[!t]
\centering
\includegraphics[width=\textwidth,trim={0 0cm 0 1.5cm},clip]{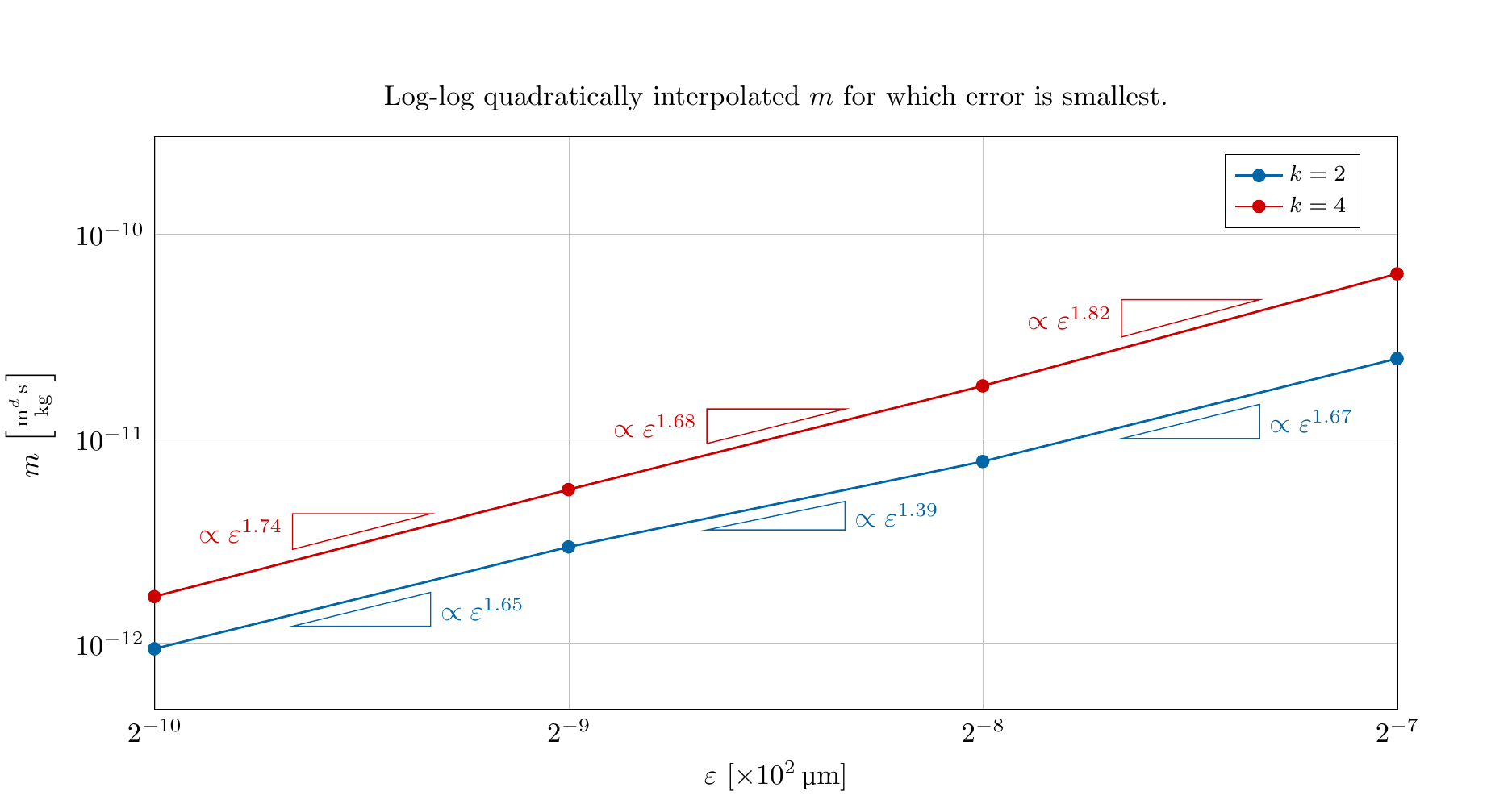}
\caption{Optimal mobility $m$ obtained from quadratic interpolation around the minima in Tables~\ref{tab:mobsk2}~and~\ref{tab:mobsk4}.}
\label{fig:m_eps_k2_k4}
\end{figure}

To provide a more precise assessment of the optimal scaling relation $m:=m_{\varepsilon}$, we determine for each~$\varepsilon$ the optimal value of~$m$ based on a quadratic log-log interpolation around the minimal values in Tables~\ref{tab:mobsk2} and~\ref{tab:mobsk4}. Figure~\ref{fig:m_eps_k2_k4} plots the optimal value of~$m$ versus~$\varepsilon$. For both wave numbers, we observe an optimal scaling~$m\propto\varepsilon^{a_{\OPT}}$ with $a_{\OPT}\approx1.7$. It is noteworthy that the constant of proportionality in the scaling relation is different for the two modes, and that the graphs are offset in the $\varepsilon$\nobreakdash-dependence by a factor of approximately two, i.e.\ the optimal mobility for~$k=4$ is approximately $2^{a_{\OPT}}$ larger than the optimal mobility for~$k=2$. This suggests that the optimal mobility in fact scales with~$(\varepsilon/\ell)^{a_{\OPT}}$, where $\ell$ represents another characteristic length scale of the interface which, for the considered droplet-oscillation test case, is proportional to the wave length of the perturbation. This observation calls for further investigation, but we consider a detailed analysis of this aspect beyond the scope of the present work.

\subsection{Sensitivity to the proportionality constant}
\label{ssec:numexperror}

In the previous section, we established optimal values of the mobility parameter and inferred an optimal scaling $m\propto\varepsilon^{a_{\OPT}}$ in the sharp-interface limit. The results in Section~\ref{ssec:numexpmobility} also convey that the constant of proportionality in the scaling relation $m_{\varepsilon}=\mathscr{C}\varepsilon^{a_{\OPT}}$ depends on the configuration and dynamics of the interface. This raises the question how sensitive the solution is to suboptimality of the proportionality constant in the scaling relation. 

\begin{figure}[!t]
\centering
\includegraphics[width=\textwidth]{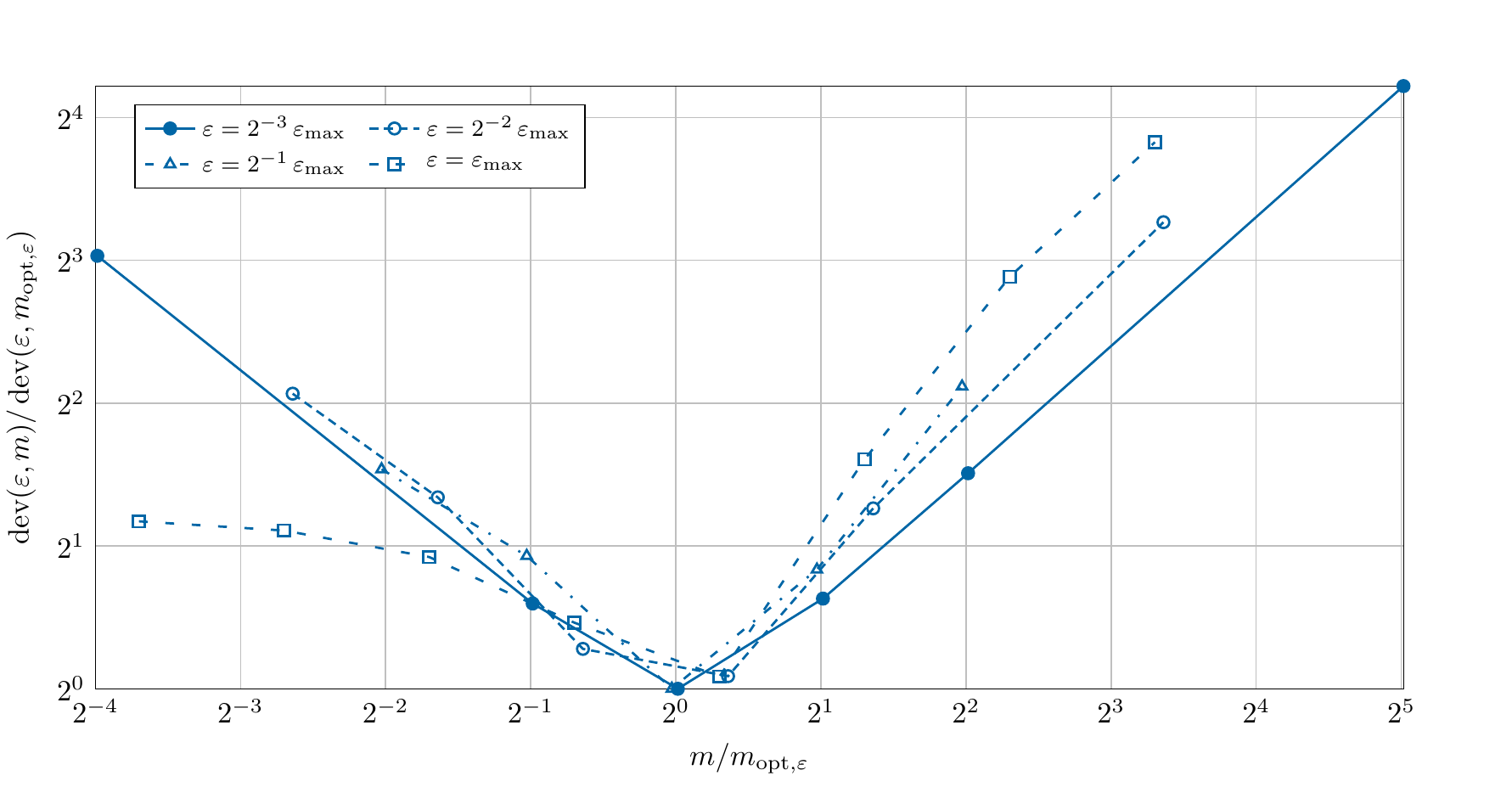}
\caption{%
Normalized deviation $\operatorname{dev}(\varepsilon,m)/\operatorname{dev}(\varepsilon,m_{\OPT,\varepsilon})$ versus normalized mobility $m/m_{\OPT,\varepsilon}$ for $k=2$.}
\label{fig:dev_m_k2}
\end{figure}

\begin{figure}[!t]
\centering
\includegraphics[width=\textwidth]{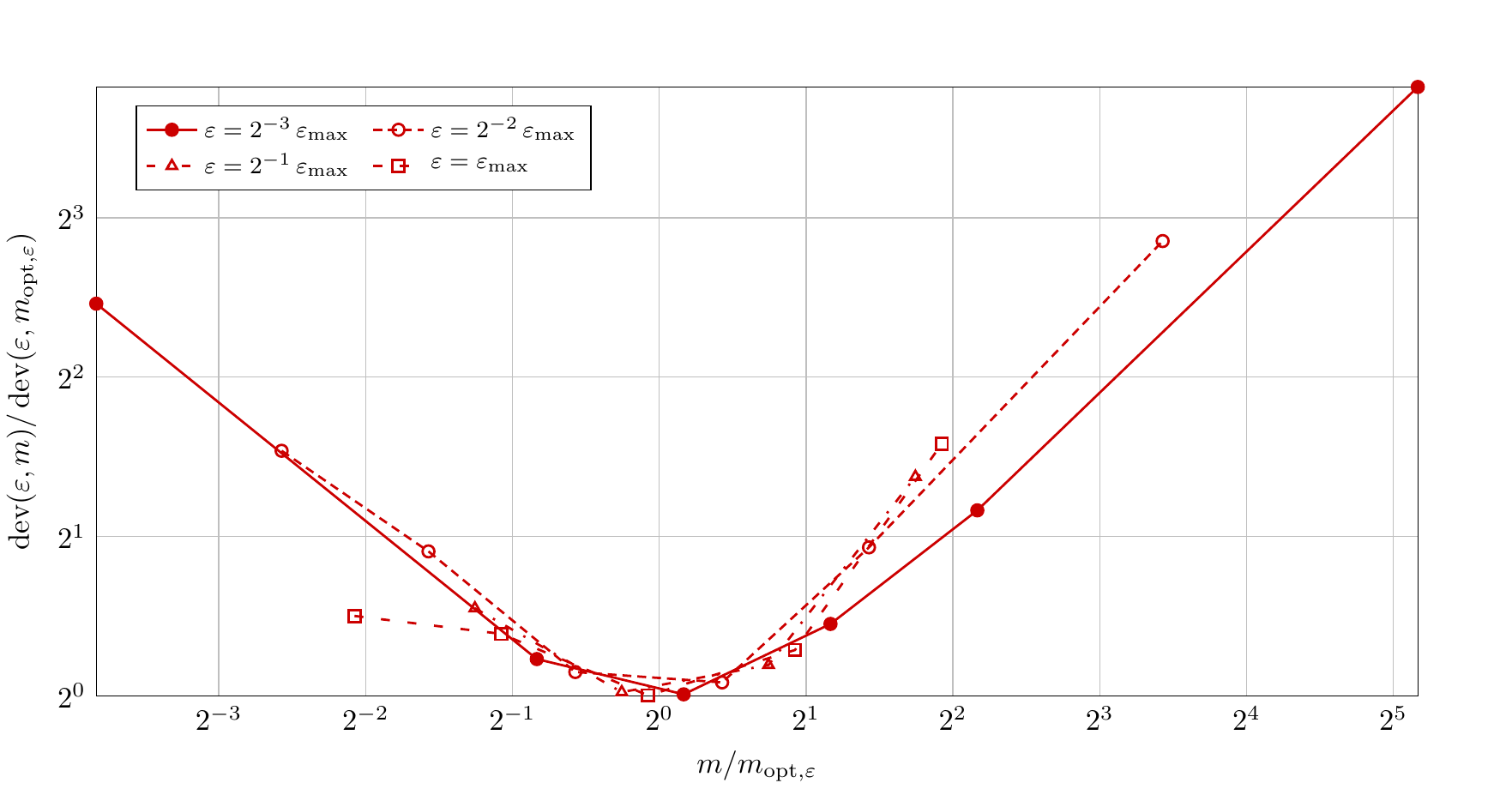}
\caption{%
Normalized deviation $\operatorname{dev}(\varepsilon,m)/\operatorname{dev}(\varepsilon,m_{\OPT,\varepsilon})$ versus normalized mobility $m/m_{\OPT,\varepsilon}$ for $k=4$.}
\label{fig:dev_m_k4}
\end{figure}

To elucidate the sensitivity of the deviation of the diffuse-interface solution to the sharp-interface solution with respect to the mobility in the limit $\varepsilon\to{}+0$, Figure~\ref{fig:dev_m_k2} (resp. Figure~\ref{fig:dev_m_k4}) plots for each $\varepsilon\in\mathcal{E}$ the ratio of the deviation 
$\operatorname{dev}(\varepsilon,m)$ in the columns of Table~\ref{tab:mobsk2} (resp. Table~\ref{tab:mobsk4})  to the minimal deviation $\operatorname{dev}(\varepsilon,m_{\OPT,\varepsilon})$ versus the ratio~$m/m_{\OPT,\varepsilon}$. Noting that the curves in Figures~\ref{fig:dev_m_k2} and~\ref{fig:dev_m_k4} exhibit a vanishing slope near $m/m_{\OPT,\varepsilon}=1$, one can conclude
that in the vicinity of the optimal mobility, the relative deviation is essentially independent of the mobility. However, for larger departures from the optimal mobility and sufficiently small~$\varepsilon$, the relative deviation increases almost linearly in~$\max(m/m_{\OPT,\varepsilon},m_{\OPT,\varepsilon}/m)$. For the largest~$\varepsilon\in\mathcal{E}$, the relative deviation appears to be less sensitive to underestimation than to overestimation of the mobility. However, the results plotted with solid markers in Figures~\ref{fig:dev_m_k2} and~\ref{fig:dev_m_k4} indicate that in the sharp-interface limit, the relative deviation is equally sensitive to under- and overestimation of the mobility.

\subsection{Suboptimal mobility scaling and convergence in the sharp-interface limit}

\begin{figure}[!b]
\centering
\includegraphics[width=\textwidth,trim={0 0cm 0 0.5cm},clip]{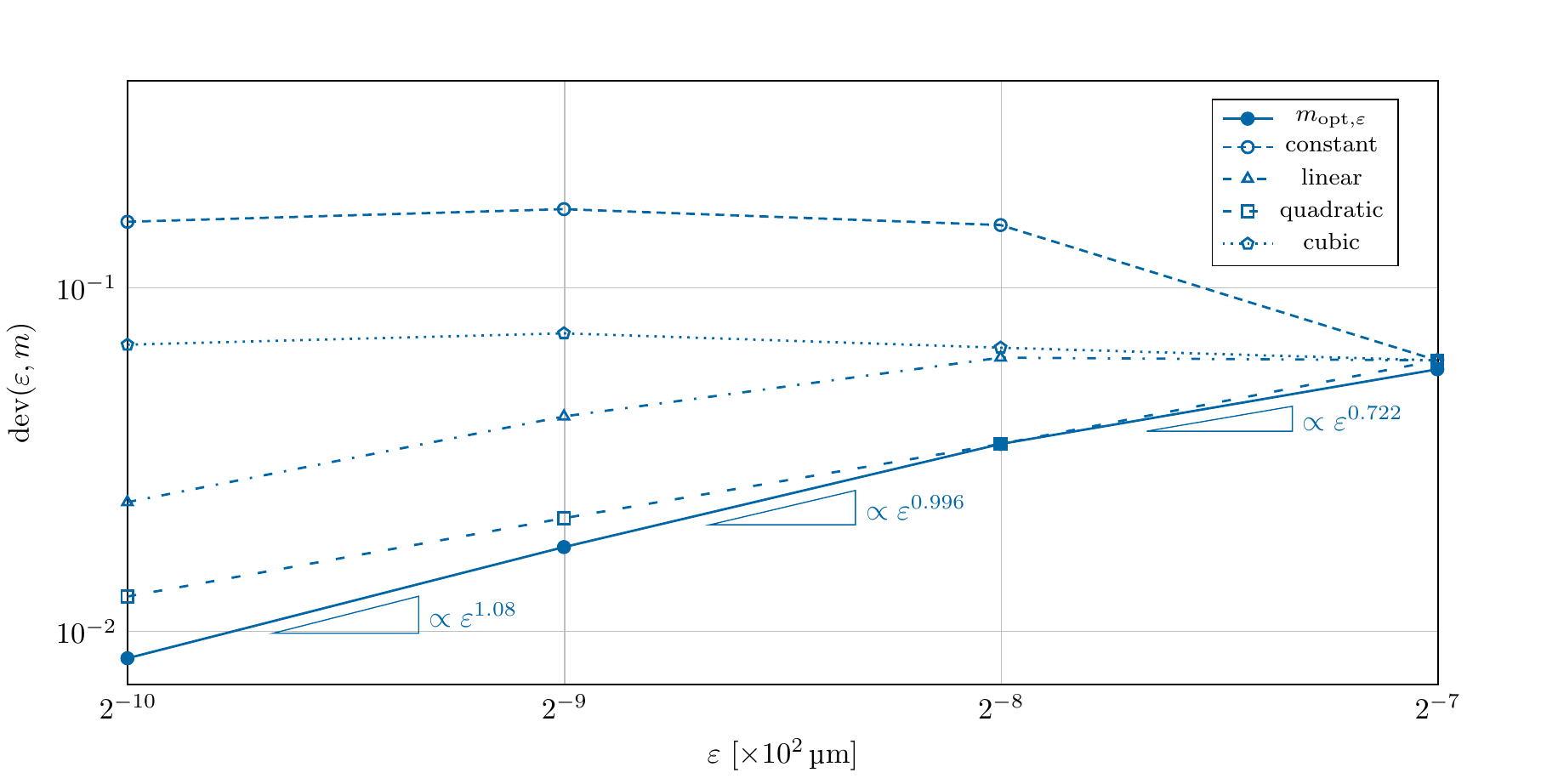}
\caption{Error convergence rates for optimal and suboptimal scalings of mobility $m$, for mode of oscillation $k=2$.}
\label{fig:dev_eps_k2}
\end{figure}

\begin{figure}[!b]
\centering
\includegraphics[width=\textwidth,trim={0 0cm 0 0.5cm},clip]{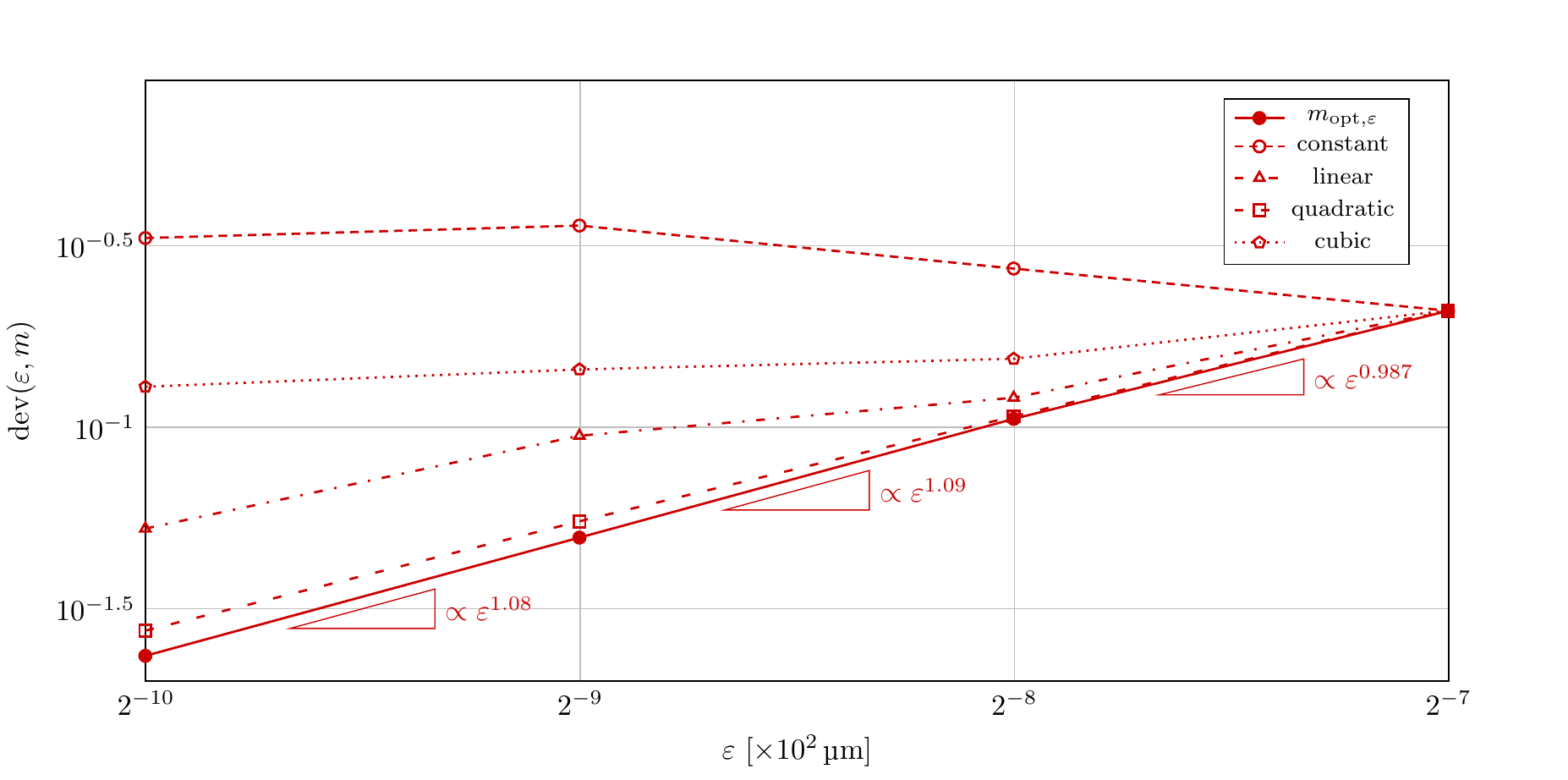}
\caption{Error convergence rates for optimal and suboptimal scalings of mobility $m$, for mode of oscillation $k=4$.}
\label{fig:dev_eps_k4}
\end{figure}

To illustrate the effect of the scaling of the mobility on the approach to the sharp-interface limit solution, Figures~\ref{fig:dev_eps_k2} and~\ref{fig:dev_eps_k4}
plot the deviation $\operatorname{dev}(\varepsilon,(\varepsilon/\varepsilon_{\MAX})^a \, m_0)$ versus~$\varepsilon$ for $a\in\{0,\ldots,3\}$. Herein, $m_0$ corresponds to the optimal sampled mobility for 
$\varepsilon_{\MAX}$, viz.\ $m_0=2^{-3}m_{\MAX}$ for $k=2$ and $m_0=2^{-2}m_{\MAX}$ for $k=4$; cf.\ Tables~\ref{tab:mobsk2} and~\ref{tab:mobsk4}. For reference, the figures also contain the estimated minimal deviation obtained by minimization of the quadratic interpolation, $\operatorname{dev}(\varepsilon,m_{\OPT,\varepsilon})$. The figures convey that for the optimal scaling of the mobility, the diffuse-interface solution approaches the sharp-interface solution essentially at order~$\varepsilon$, i.e.\ $\operatorname{dev}(\varepsilon,m_{\OPT,\varepsilon})=\mathcal{O}(\varepsilon)$ as $\varepsilon\to{}+0$. For the linear and quadratic scaling of the mobility with the interface thickness, $m\propto\varepsilon^a$ with $a\in\{1,2\}$, i.e.\ the integer scalings of the mobility closest to the optimum, we observe convergence to the sharp-interface solution, but at a suboptimal (sublinear) rate. For the constant and cubic scalings of the mobility, $m \propto \varepsilon^a$ with $a\in\{0,3\}$, the deviation $\operatorname{dev}(\varepsilon,m_0(\varepsilon/\varepsilon_{\MAX})^a)$ does not vanish as $\varepsilon\to{}+0$, i.e.\ the diffuse-interface solution does not convergence to the sharp-interface solution. For $m \propto \varepsilon^0$, this confirms the known result that for constant mobility, the Abels--Garcke--Gr\"{u}n NSCH model converges to the nonclassical sharp-interface Navier--Stokes/Mullins--Sekerka model; see~\cite{Abels:2018ly}.

\begin{remark}
It is noteworthy that the subcubic scaling of the mobility, $m\propto\varepsilon^a$ with $0<{}a<3$, which is necessary to converge to the classical sharp-interface solution in the limit \mbox{$\varepsilon\to{}+0$}, implies that the characteristic diffusive time scale $T_{\mathrm{diff}}:=\varepsilon^3/\sigma{}m$ associated with the diffuse interface, approaches zero in the sharp-interface limit. Consequently, for any $\varepsilon$\nobreakdash-independent characteristic time scale~$T_*$ in the problem under consideration, e.g.\ the period of oscillation of a droplet, it holds that
the ratio $T_{\mathrm{diff}}/T_*\to{}+0$ as $\varepsilon\to{}+0$. For numerical time-integration methods for the NSCH equations, it is therefore essential that such methods are robust in the limit 
$T_{\mathrm{diff}}/\tau\to+0$ (with~$\tau$ denoting the time-step size), to avoid excessive computational complexity in the sharp-interface limit.  
\end{remark}

%===============================================================================================================%
%===============================================================================================================%
%===============================================================================================================%

\section{Conclusions}\label{sec:concl}
Diffuse-interface binary-fluid models bear significant potential for describing complex phenomena in fluid mechanics, such as topological changes of the fluid-fluid interface and dynamic wetting, by virtue of their implicit representation of the interface. In the absence of topological changes of the interface, diffuse-interface models should reduce to corresponding classical sharp-interface models in the so-called sharp-interface limit, viz.\ if the interface-thickness parameter,~$\varepsilon$, passes to zero. Contemporary understanding of the sharp-interface limit is however incomplete and, in particular, the scaling of the mobility parameter, $m$, with~$\varepsilon$ as $\varepsilon\to{}+0$ is incompletely understood. In this article, we investigated the limit behavior of the Abels--Garcke--Gr\"un Navier-Stokes--Cahn-Hilliard model for the classical case of an oscillating droplet in two dimensions, by means of an adaptive finite-element methodology.

To provide reference sharp-interface solutions, we derived new two-dimensional analytical expressions for the velocity and pressure fields for small-amplitude oscillations of a viscous droplet in a viscous ambient fluid with different densities and viscosities.

For mode numbers $k=2,4$, we compared the solutions of the Navier-Stokes--Cahn-Hilliard model to the corresponding analytical solutions for a decreasing sequence of interface-thickness parameters and a suitably chosen sequence of mobility parameters. Based on an analysis of the deviation between the diffuse-interface solution and the sharp-interface solution, we deduced that $m_{\OPT,\varepsilon}\propto\varepsilon^{a_{\OPT}}$ with ${a_{\OPT}}\approx1.7$ corresponds to the optimal scaling of the mobility in the sharp-interface limit. We found that this optimal scaling is universal for $k=2$ and $k=4$. However, we also observed that the factor of proportionality differs by a factor of approximately $2^{a_{\OPT}}$, indicating that the optimal mobility in fact scales with the ratio of~$\varepsilon$ to another characteristic length scale, in this case, proportional to the wave length of the perturbation. The observed dependence of the optimal mobility on the configuration and motion of the interface warrants further investigation.

For the optimal scaling of the mobility parameter, we observed that the deviation between the diffuse-interface solution and its sharp-interface limit decreases according to~$\mathcal{O}(\varepsilon)$ in the limit $\varepsilon\to{}+0$. Our investigation of suboptimal integer scalings of the mobility conveyed that the sharp-interface limit is also attained for a linear and quadratic scaling of the mobility, but not for a constant or cubic scaling. For the linear and quadratic scaling, the approach of the diffuse-interface solution to the sharp-interface solution occurs at a suboptimal (sublinear) rate. The fact that the scaling of the mobility with~$\varepsilon$ must be subcubic, implies that the characteristic diffusive time scale $\varepsilon^3/\sigma{}m$ (with~$\sigma$ as surface tension) passes to zero in the sharp-interface limit.

\section*{Acknowledgments}
This research was partly conducted within the Industrial Partnership Program {\it Fundamental Fluid Dynamics Challenges in 
Inkjet Printing\/} ({\it FIP\/}), a joint research program of Canon Production Printing, Eindhoven University of Technology,
University of Twente, and the Netherlands Organization for Scientific Research (NWO). T.H.B.\ Demont and S.K.F.\ Stoter gratefully acknowledge financial support through the FIP program. All simulations have been performed using the open source software package Nutils  (www.nutils.org) \cite{nutils}.

The authors thank Dr.\ S.W.\ Rienstra of the Department of Mathematics and Computer Science at Eindhoven University of Technology for insightful discussions with respect to the interpretation of matched-asymptotic-expansion results for sharp-interface limits of the NSCH equations in the literature.
\appendix

\bibliographystyle{plain}
\bibliography{osc_drop}

\end{document}